\def\ACKS{}
\def\mms{}

\documentclass[12pt]{article}
\usepackage{natbib}
\usepackage{url}
\usepackage{geometry}
\usepackage{amsmath,amsthm,amsfonts,breqn}

\usepackage{graphicx}
\graphicspath{{./graphics/}}
\usepackage[shortcuts]{extdash}
\usepackage{csvsimple}

\usepackage {tikz}
\usetikzlibrary {positioning}

\usepackage{enumitem}
\setlist{nolistsep}

\theoremstyle{plain}
\newtheorem{maintheorem}{Theorem}
\newtheorem{corollary}{Corollary}

\newtheorem{theorem}{Theorem}[section]
\newtheorem{lemma}[theorem]{Lemma}
\newtheorem{claim}[theorem]{Claim}
\newtheorem{proposition}[theorem]{Proposition}
\newtheorem*{efmtheorem*}{EFM theorem}
\newtheorem*{theorem*}{Theorem}
\newtheorem{open}{Open question}

\theoremstyle{definition}
\newtheorem{remark}[theorem]{Remark}
\newtheorem{definition}[theorem]{Definition}

\usepackage{algorithm}
\usepackage{algorithmic}

\newcounter{problem}
\makeatletter
\newenvironment{problem}[1][htb]{%
    \let\c@algorithm\c@problem
    \renewcommand{\ALG@name}{Problem}
    \begin{algorithm}[#1]%
    }{\end{algorithm}
}
\makeatother

\newcommand{\range}[2]{\in\{#1,\dots,#2\}}
\newcommand{\bestk}[1]{\textrm{Best}^k_{#1}} 
\newcommand{\mmsk}[1]{\textrm{MMS}^{k,n}_{#1}} \newcommand{\mmskm}[1]{\textrm{MMS}^{k,{n-1}}_{#1}} 
\newcommand{\vbestk}[2]{V_{#1}\left(\bestk{#1}(#2)\right)}
\newcommand{\dottedline}{\hbox to 15cm{\leaders\hbox to 5pt{\hss.\hss}\hfil}}

\def\prop(#1,#2,#3){\operatorname{\textsc{Prop}}(#1,#3,#2)}
\def\prope(#1,#2,#3){\operatorname{\textsc{PropEF}}(#1,#3,#2)}
\def\proprel(#1,#2,#3){\operatorname{\textsc{RelProp}}(#1,#3,#2)}
\def\properel(#1,#2,#3){\operatorname{\textsc{RelPropEF}}(#1,#3,#2)}

\title{Fair Multi-Cake Cutting}

\begin{document}

\author{Erel Segal-Halevi
\\
Ariel University
\\
Kiriat Hamada 3, Ariel 40700, Israel
\\
erelsgl@gmail.com
}
\date{}
\maketitle

\begin{abstract}
In the classic problem of fair cake-cutting, a single interval (``cake'') has to be divided among $n$ agents with different value measures, giving each agent a single sub-interval with a value of at least $1/n$ of the total.
This paper studies a generalization in which the cake is made of m disjoint intervals, and each agent should get at most k sub-intervals.
The paper presents a polynomial-time algorithm that guarantees to each agent at least $\min ( 1/n, k/(m+n-1) )$ of the total value, and shows that this is the largest fraction that can be guaranteed.
The algorithm simultaneously guarantees to each agent at least $1/n$ of the value that the agent can get without partners to share with.
The main technical tool is envy-free matching in a bipartite graph.
Some of the results remain valid even with additional fairness constraints such as envy-freeness.

Besides the natural application of the algorithm to simultaneous division of multiple land-estates, the paper shows an application to a geometric problem --- fair division of a two-dimensional land estate shaped as a rectilinear polygon, where each agent should receive a rectangular piece.
\end{abstract}


\newpage
\section{Introduction}
\label{sec:intro}
Consider $n$ people who inherit $m$ distant land-estates (henceforth ``islands'') and want to divide the property fairly among them. 
Classic algorithms for \emph{proportional cake-cutting} \citep{Steinhaus1948Problem,Even1984Note}
can be used to divide each island into $n$ pieces such that the value of piece $i$, according to the personal value-measure of person $i$, is at least $1/n$ of the total island value.
However, this scheme requires each person to manage properties on $m$ distant locations, which may be quite inconvenient. 
An alternative scheme is to consider the union of all $m$ islands as a single cake (henceforth ``multicake''), and partition it into $n$ pieces using the above-mentioned  algorithms. However, this too might give some agents a share that overlaps many distinct islands. For example, with $m=5$ islands and $n=3$ agents, a typical partition might look like this:

\begin{center}
\includegraphics[width=0.7\textheight]{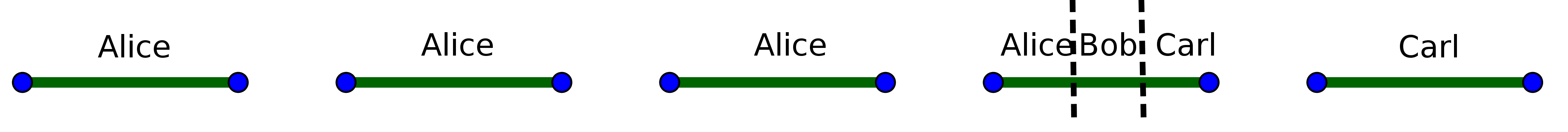}
\end{center}
where Alice's share contains 4 disconnected pieces and Carl's share contains 2 disconnected pieces.%
\footnote{
Such a partition may occur if the three leftmost islands have a low value while the fourth island is very valuable.
}
This may require Alice to waste a lot of time flying between her different estates.

Can we find a more convenient division? For example, 
can we divide the multicake fairly such that each agent receives at most $3$ disconnected pieces? In general, what is the smallest number $k$ (as a function of $m$ and $n$) such that there always exists a fair division of the multicake in which each agent receives at most $k$ disconnected pieces?

This paper studies the following more general question. We are given an integer parameter $k\geq 1$, and it is required to give each agent a share made of at most $k$ disconnected pieces. What is the largest fraction $r(m,k,n)$ such that there always exists an allocation giving each agent at least a fraction $r(m,k,n)$ of his/her total property value?

An obvious upper bound is $r(m,k,n)\leq 1/n$. It is attainable when $k\geq m$, since then we can just divide each island separately using classic cake-cutting algorithms.
The question is more challenging when $k<m$.
The main result of this paper is:
\begin{maintheorem}
\label{thm:prop}
For every $m\geq 1, k\geq 1, n\geq 1$:\hskip 1cm
$
r(m,k,n) = \min\left({1\over n},{k\over m+n-1}\right).
$
\end{maintheorem}
In words: it is always possible to guarantee to each agent at least a fraction 
$\min\left({1\over n},{k\over m+n-1}\right)$ of the total multicake value, and in some cases it is impossible to guarantee a higher fraction.

The upper bound $r(m,k,n) \leq \min\left({1\over n},{k\over m+n-1}\right)$ is proved by a simple example (\textbf{Section \ref{sec:negative}}).
The matching lower bound is proved constructively  by a polynomial-time algorithm, which is the main technical contribution of the paper (\textbf{Section \ref{sec:algorithm}}). 
The case $k=1$ is simple and can be handled using classic cake-cutting techniques (sub. \ref{sub:k=1}). However, for $k\geq 2$ these techniques yield a sub-optimal result, and to match the upper bound we need a new technique using an \emph{envy-free matching in a bipartite graph} (sub. \ref{sub:k=2}).

Whenever $k\geq 1 + {m-1 \over n}$,
Theorem \ref{thm:prop} implies that $r(m,k,n)=1/n$. Therefore,
\begin{corollary}
If $k\geq 1 + {m-1\over n}$ (each agent is willing to accept at least $1 + {m-1\over n}$ disconnected pieces), then a proportional multicake allocation exists and can be computed efficiently.
\end{corollary}
In particular, in the opening example, it is indeed possible to find a fair allocation in which each person receives at most $3$ disconnected pieces.

Theorem \ref{thm:prop} provides an \emph{absolute} value guarantee --- a fraction of the total multicake value. 
One may also be interested in a \emph{relative} value guarantee --- a fraction of the value that an agent could attain without partners, which is the value of his/her $k$ most valuable islands. 
What is the largest fraction $r_R(m,k,n)$ such that there always exists an allocation giving each agent at least a fraction $r_R(m,k,n)$ of his/her $k$ best islands?
The next result shows that this fraction is $1/n$.%
\begin{maintheorem}
\label{thm:relprop}
For every $m\geq 1, k\geq 1, n\geq 1$:
\hskip 1cm 
$
r_R(m,k,n) = 1/n.
$
\end{maintheorem}
The upper bound of $1/n$ is obvious.%
\footnote{
For example, consider $n$ agents with the same valuation, all of whom value a single island at $1$ and the other $m-1$ islands at $0$. 
}
The lower bound is proved using the same algorithm as Theorem \ref{thm:prop}.
Moreover, the algorithm can even allow each agent to choose between the two guarantees. This may be useful, as for each agent, a different guarantee may be preferred.
As an example, suppose there are $n=10$ agents and $m=11$ islands and $k=1$.
Theorem \ref{thm:prop} guarantees to each agent a fraction $1/20$ of the total value.
Suppose Alice values all $11$ islands uniformly at $1/11$ of the total value, 
while George thinks that all the multicake value is concentrated at a single island.
Theorem \ref{thm:relprop} guarantees Alice 
$1/110$ and George $1/10$ of the total value.
Our algorithm can give the better guarantee to each agent, i.e., at least
 $1/20$ to Alice and at least $1/10$ to George (\textbf{Section \ref{sec:relprop}}).

Besides proportionality, 
an additional fairness condition is \emph{envy-freeness}. It requires that every agent values his/her share at least as much as the share of every other agent.
In the classic setting where $m=k=1$, 
it is known that envy-freeness implies proportionality,
and an envy-free (and proportional) allocation always exists \citep{Woodall1980Dividing,Stromquist1980How,Su1999Rental}.
When $n=2$, an envy-free allocation can be computed efficiently by the classic cut-and-choose protocol; 
when $n\geq 3$, an envy-free allocation cannot be computed by a finite protocol \citep{Stromquist2008Envyfree}.
What is the largest fraction $r_E(m,k,n)$ such that there always exists an \emph{envy-free} allocation giving each agent at least a fraction $r_E(m,k,n)$ of his/her total property value?

An obvious upper bound is $r_E(m,k,n)\leq r(m,k,n)$: adding a constraint on the allocations cannot help  attain a higher value-guarantee.
A lower bound is presented below.
\begin{maintheorem}
\label{thm:envyfree}
For every $m\geq 1, k\geq 1$:

(a) When $n=2$:
\hskip 1cm
$r_E(m,k,n) = \min\left({1\over n},{k\over m+n-1}\right).$

(b) For any $n\geq 3$:
\hskip 1cm
$r_E(m,k,n) \geq 
\min\left({1\over n},{k\over m+nk-k}\right).$
\end{maintheorem}
Part (a) is proved constructively by a polynomial-time algorithm.
It exactly matches the upper bound of Theorem \ref{thm:prop}: when there are two agents, 
adding an envy-freeness constraint does not decrease the value-guarantee.

Part (b) is proved  
existentially by a reduction to the classic setting (\textbf{Section \ref{sec:envyfree}}). 
The bound is tight for $k=1$ but not for $k\geq 2$.
The exact value of $r_E(m,k,n)$ 
for $n\geq 3$ and $k\geq 2$ remains unknown.

Besides the natural application of the multicake division algorithm to simultaneous division of multiple land-estates, it can be applied to 
the problem of dividing a two-dimensional land estate shaped as a rectilinear polygon, when each agent requires a rectangular land-plot. If the land-estate has $T$ reflex vertices (vertices with internal angle $270^{\circ}$), then each agent can be allocated a rectangle worth at least $1/(n+T)$ of the total polygon value, and this is tight (\textbf{Section \ref{sec:rectilinear}}).

Various ideas for future work are presented at ending subsections of sections 
 \ref{sec:algorithm},
\ref{sec:relprop},
 \ref{sec:envyfree},
  \ref{sec:rectilinear}.

\subsection{Related work}
Fair division problems have inspired hundreds of research papers and several books \citep{Brams1996Fair,Robertson1998CakeCutting,Moulin2004Fair,Brams2007Mathematics}. 
See 
\citet{Procaccia2015Cake,Branzei2015Computational,Lindner2016,Brams2017FairDivision,SegalHalevi2017Phd}
for recent surveys.%
\ifdefined\ACKS
\footnote{
Some relatively up-to-date surveys can be found in the Wikipedia pages ``Fair division'', ``Fair cake-cutting'', ``Proportional cake-cutting'', ``Envy-free cake-cutting'' and ``Fair item allocation''.}
\fi
~Some works that are more closely related to the present paper are surveyed below.

Most works on cake-cutting either require to give each agent a single connected piece, or ignore connectivity altogether; few works consider the natural intermediate case, in which each agent should receive at most a fixed number $k$ of disconnected pieces.
We are aware of two exceptions.

(a)
\citet{ram2018price} consider a connected cake, i.e, a single island. They quantify the gain, as a function of $k$, in the optimal \emph{social welfare}, defined as the sum of values of all agents or the minimum value per agent. Based on their results, they conjecture that this gain grows linearly with $\min(n,k)$. 

(b)
\citet{bei2019dividing} consider a cake modeled as a connected graph, where the agents prefer to receive a connected subset of edges. This is a generalization of the classic model of an interval cake. They prove that, for every connected graph, it is possible to give each agent a connected piece worth at least ${1 \over 2n-1}$ of his/her total value. For a star with $m<2n-1$ edges, the fraction improves to $\frac{1}{n + \lceil m/2 \rceil - 1}$.
When there are $n=2$ agents, and the agents may get disconnected pieces whose total count (for both agents together) is at most $k + 1$, the fraction improves to ${1\over 2}\cdot (1 - {1\over 3^k})$.

The present paper is more general in that it considers a cake that may itself be disconnected.
It is less general in that each connected component must be an interval (i.e., a single edge).
A potential future work avenue is to generalize both papers by letting the cake be an arbitrary (possibly disconnected) graph.

Some works aim to minimize the number of cuts required to attain various fairness and efficiency goals 
\citep{Webb1997How,Shishido1999MarkChooseCut,Barbanel2004Cake,Barbanel2014TwoPerson,alijani2017envy,seddighin2019expand,SegalHalevi2018Entitlements}. These works usually assume that the cake itself is connected. Moreover, their goal is to minimize the global number of cuts, while the goal in the present paper is to ensure that every individual agent is not given too many disconnected crumbs.
Nevertheless, our algorithm makes at most $n-1$ cuts, which in some cases is the smallest possible amount (see Remark \ref{rem:cuts}).

Some works consider a different multi-cake division problem, in which each agent must get a part of \emph{every} sub-cake \citep{Cloutier2010Twoplayer,Lebert2013Envyfree,nyman2020fair}. This is opposite to the present paper, in which each agent wants to overlap as few cakes (islands) as possible.

Some works consider an endowment made of several discrete objects, and try to minimize the number of items that have to be shared.
\citet{brams2000winwin} present the Adjusted Winner procedure, which divides items fairly and efficiently among two agents such that at most one item is shared. 
\citet{wilson1998fair}, in an unpublished manuscript, shows that this can be generalized to $n$ agents such that at most $n-1$ items are shared, and this is the smallest number that can be guaranteed. 
\citet{segal2019fair} extends these bounded sharing results to other fairness criteria, while 
\citet{sandomirskiy2019fair} present an algorithm for finding a fair and efficient division with minimal sharing.
The present paper focuses on the convenience  of each individual agent, rather than the global number of items (islands) that are shared.

\section{Preliminaries}
There is a set $C$ (``the multicake'') that contains $m\geq 1$ disjoint subsets (``the islands'').
It is convenient to assume that $C$ is a subset of the real line and the islands are pairwise-disjoint intervals of unit length.
Given such an island $B\subseteq C$ and a real number $d\in[0,1]$, we denote by $B[d]$ the sub-interval of length $d$ at the left-hand side of $B$.

There are $n\geq 1$ agents. 
The preferences of each agent $i$ are represented by a non-negative integrable function $v_i: C\to\mathbb{R}_+$ called a \emph{value-density}.
The agent's value of each piece $Z_i \subseteq C$ 
is denoted by $V_i(Z_i)$, which is defined as the integral of the agent's value-density: 
\begin{align*}
V_i(Z_i) := \int_{x\in Z_i} v_i(x) dx.
\end{align*}
The definition implies that each $V_i$ is a \emph{non-negative measure} (an additive set function) on $C$, and that it is \emph{non-atomic} --- the value of zero-length subsets is zero.

An \emph{allocation} $\mathbf{Z} = Z_1,\ldots,Z_n$ is a vector of $n$ pairwise-disjoint subsets of $C$. 

There is a fixed integer constant $k\geq 1$, which denotes the maximum number of disconnected pieces that an agent can use. 
An allocation $\mathbf{Z}$ is called \emph{$k$-feasible} if for each $i\in[n]$, $Z_i$ is the union of at most $k$ intervals. In particular, $Z_i$ overlaps at most $k$ different islands.

Free disposal is assumed, i.e., some parts of $C$ may remain unallocated. This assumption is unavoidable even with $n=1$ agent. For example, if there are $m=2$ islands, then an agent who wants $k=1$ interval must leave one island unallocated.

For every subset $Z\subseteq C$, $\bestk{i}(Z)$ denotes the $k$ intervals contained in $Z$ that are most valuable for agent $i$ (breaking ties arbitrarily).
In particular, $\bestk{i}(C)$ denotes the $k$ islands most valuable for $i$.
If an allocation $\mathbf{Z}$ is $k$-feasible,
then every subset $Z_i\subseteq C$ is made of at most $k$ intervals, so $\bestk{i}(Z_i) = Z_i$. In general $V_i(\bestk{i}(Z_i))$ may be smaller than $V_i(Z_i)$.

\section{Upper Bound}
\label{sec:negative}
Before presenting an algorithm for dividing a  multicake, it is useful to have an upper bound on what such an algorithm can hope to achieve.
\begin{lemma}
\label{lem:upperbound}
For every $m\geq 1, n\geq 1, k\geq 1$, 
there is a multicake instance with $m$ islands and $n$ agents in which, in every $k$-feasible allocation, at least one agent receives at most the following fraction of the total multicake value:
\begin{align*}
\min\left({1\over n},{k \over  n+m-1}\right)
\end{align*}
\end{lemma}
\begin{proof}
Consider a multicake with $m-1$ ``small'' islands and one ``big'' island.
All $n$ agents have the same value-measure: they value every small island at $1$ and value the big island at $n$. So the total multicake value for every agent $i$ is  $V_i(C) = n+m-1$. 

Since the valuations are identical, it is obviously impossible to give every agent $i$ a value of more than $V_i(C)/n$.
It remains to show that it is impossible to give every agent a value of more than $k$.
~~
Consider the following two cases:
\begin{enumerate}
\item At least one agent gets all his/her $k$ pieces in some $k$ small islands. Since the value of every small island is 1, this agent receives a value of at most $k$.
\item All $n$ agents get at least one of their $k$ intervals in the big island. 
Since the value of the big island is $n$,
at least one agent receives a value of at most $1$ from that island. 
This agent receives a value of at most $k-1$ from his/her other $k-1$ pieces, which are subsets of small islands. Therefore his/her total value is at most $k$.
\end{enumerate}
In both cases, at least one agent receives a value of at most ${k\over n+m-1}V_i(C)$.
\end{proof}

\section{Algorithm}
\label{sec:algorithm}
Motivated by the upper bound of Lemma \ref{lem:upperbound}, this section aims to solve the following problem.

\begin{problem}
\caption{
\label{prob:multicake}
Multicake Division Problem 
}
\begin{algorithmic}
\REQUIRE ~\\

\begin{itemize}[noitemsep]
\item Positive integers $m, n, k$; 
\item A multicake $C$ with $m$ islands
\item $n$ value-measures $V_1,\ldots,V_n$ on $C$.
\end{itemize}
\ENSURE A $k$-feasible allocation $\mathbf{Z} = Z_1,\ldots,Z_n$ in which for all $i\in[n]$:
\begin{align*}
{V_i(Z_i)} \geq \min\left({1\over n},{k \over  n+m-1}\right) \cdot  V_i(C).
\end{align*}
\end{algorithmic}
\end{problem}

Below, \S\ref{sub:main} presents the main algorithm: it is a recursive algorithm that, in each iteration, finds a partial allocation of a subset of $C$ to a subset of the agents.
The main challenge is to find the appropriate partial allocation in each step.
As a warm-up, \S\ref{sub:k=1} and \S\ref{sub:k=2} show how to find a partial allocation in the cases $k=1$ and $k=2$ respectively. 
\S\ref{sub:k=k} presents the algorithm for finding a partial allocation for any $k$, thus completing the description of the multicake division algorithm.
\S\ref{sub:example} presents an example, which the reader is invited to follow in parallel to the algorithm.
\S\ref{sub:algorithm-open} presents some open questions.

\subsection{Main loop}
\label{sub:main}
We first reduce Problem \ref{prob:multicake} into a more convenient normalized form displayed below as Problem \ref{prob:multicake-normalized}.

\begin{problem}
\caption{
\label{prob:multicake-normalized}
Normalized Multicake Division Problem 
}
\begin{algorithmic}
\REQUIRE ~\\
\begin{itemize}[noitemsep]
\item Positive integers $m, n, k$ with $m\geq nk-n+1$;
\item A multicake $C$ with $m$ islands; 
\item $n$ value-measures $V_1,\ldots,V_n$ on $C$ having $V_i(C) \geq n+m-1$ for all $i\in[n]$.
\end{itemize}
\ENSURE A $k$-feasible allocation $\mathbf{Z} = Z_1,\ldots,Z_n$ in which for all $i\in[n]$:
\begin{align*}
V_i(Z_i)\geq k
\end{align*}
\end{algorithmic}
\end{problem}

\begin{lemma}
\label{lem:normalization}
Any instance of problem \ref{prob:multicake} can be reduced to an instance of problem \ref{prob:multicake-normalized}.
\end{lemma}
\begin{proof}
Given an instance of problem \ref{prob:multicake}, construct an instance of problem \ref{prob:multicake-normalized} with the same $n,k$ and with $m' = \max(m, nk-n+1)$;
if $m<nk-n+1$, then add dummy islands whose value for all agents is $0$, such that the total number of islands becomes $nk-n+1$.

For each $i\in[n]$, construct a normalized value function $V_i'(\cdot) := \dfrac{n+m'-1}{V_i(C)}\cdot V_i(\cdot)$,
so that $V_i'(C) = n+m'-1$.
Solve Problem \ref{prob:multicake-normalized} with integers $m', n, k$ and value-measures $V_1',\ldots,V_n'$; note that these inputs satisfy the input requirements of problem \ref{prob:multicake-normalized}.

The solution is an allocation $\mathbf{Z}$ such that $V_i'(Z_i)\geq k$ for all $i\in[n]$. Consider two cases:

\begin{enumerate}
\item In the original instance $m \geq nk-n+1$. Then $m' = m$, and for each agent $i\in[n]$ we have
\begin{align*}
{V_i(Z_i) \over V_i(C)} = {V_i'(Z_i) \over n+m'-1} = {V_i'(Z_i) \over n+m-1} \geq {k\over n+m-1}.
\end{align*}
\item In the original instance
$m < nk-n+1$. Then $m' = nk-n+1$, and for each agent $i\in[n]$ we have
\begin{align*}
{V_i(Z_i) \over V_i(C)} = {V_i'(Z_i) \over n+m'-1} = {V_i'(Z_i) \over nk} \geq {k\over nk} \geq {1\over n}.
\end{align*}
\end{enumerate}
In both cases, $\mathbf{Z}$ is a solution to problem \ref{prob:multicake}.
\end{proof}

Thus, to solve problem \ref{prob:multicake} it is sufficient to solve problem \ref{prob:multicake-normalized}.
We will solve problem
 \ref{prob:multicake-normalized}
by a recursive algorithm that is presented as Algorithm \ref{alg:main}.

At each recursion round, the algorithm allocates some subsets of $C$ to agents in some subset $L$. The agents in $L$ have received their due share and are therefore removed from further consideration. The remaining agents continue to divide the remainder of $C$ among them. In each round, the algorithm has to solve Problem \ref{prob:partial} below.

\begin{problem}
\caption{
\label{prob:partial}
Partial Allocation Problem.
}
\begin{algorithmic}
\REQUIRE ~\\
\begin{itemize}[noitemsep]
\item Positive integers $m, n, k$ with $m\geq nk-n+1$;
\item A multicake $C$ with $m$ islands; 
\item $n$ value-measures $V_1,\ldots,V_n$ on $C$ having $V_i(C) \geq n+m-1$ for all $i\in[n]$.
\end{itemize}
\ENSURE
A nonempty set of agents $L\subseteq [n]$, and a set of pieces $(Z_i)_{i\in L}$ with for all $i\in L$:
\begin{itemize}
\item The piece $Z_i$ contains exactly $k-1$ whole islands of $C$, plus possibly the leftmost part of some $k$-th island (so $Z_i$ contains at most $k$ intervals).
\item $V_i(Z_i)\geq k$.
\item $V_j(Z_i)\leq k$ for all $j\not\in L$.
\end{itemize}
\end{algorithmic}
\end{problem}

\begin{algorithm}[t]
\caption{
\label{alg:main}
Multicake Division: Main Algorithm
}

\begin{algorithmic}[1]
\REQUIRE ~\\
\begin{itemize}[noitemsep]
\item Positive integers $m, n, k$ with $m\geq nk-n+1$;
\\
\item A multicake $C$ with $m$ islands; 
\\
\item $n$ value-measures $V_1,\ldots,V_n$ on $C$ having $V_i(C) \geq n+m-1$ for all $i\in[n]$.
\end{itemize}
\ENSURE A $k$-feasible allocation $\mathbf{Z} = Z_1,\ldots,Z_n$ in which for all $i\in[n]$:
\begin{align*}
V_i(Z_i)\geq k
\end{align*}
\dottedline{}

\STATE 
\label{step:base}
If $n=1$, then give the single remaining agent the $k$ islands of $C$ that are most valuable to him/her, and terminate.

\STATE 
\label{step:partial}
Solve the Partial Allocation problem \ref{prob:partial}, getting a set of agents $L\subseteq [n]$ and a set of pieces $(Z_i)_{i\in L}$.
Allocate to each agent $i\in L$ the piece $Z_i$.

\STATE 
\label{step:complete}
Let $\ell := |L|$.
If 
the number of whole islands removed at step \ref{step:partial} is 
more than $\ell\cdot (k-1)$ (for example, if the $k$-th parts given to the $\ell$ agents together consume a whole island),
then add dummy islands of value $0$ such that exactly $m-\ell\cdot (k-1)$ islands remain.

\STATE 
\label{step:recurse}
Recursively partition the remainder $C\setminus (\cup_{i\in L}Z_i)$ among the remaining agents $[n]\setminus L$.
\end{algorithmic}
\end{algorithm}
            
\begin{lemma}
\label{lem:partial-implies-multicake}
Given a correct algorithm for solving the Partial Allocation problem \eqref{prob:partial}, 
the Multicake Division problem \eqref{prob:multicake-normalized}
is correctly solved by Algorithm \ref{alg:main}.
\end{lemma}
\begin{proof}
By induction on $n$. If $n=1$ then the algorithm terminates at step \ref{step:base}. The agent values $C$ at least $m+n-1 = m$, so the average value per island is at least $1$. Therefore the $k$ most valuable islands are worth for this agent at least $k$.

Suppose now that $n>1$
and that the claim holds whenever there are less than $n$ agents.
At step \ref{step:partial}, at least one agent receives a piece. Any allocated piece overlaps at most $k$ islands --- $k-1$ whole islands and possibly a piece of a $k$-th island. Hence the partial allocation is $k$-feasible.
Each allocated piece is worth at least $k$ to its receiver.

Let $\ell := |L|$ and $m' := m - \ell\cdot (k-1)$.
After $\ell$ agents receive $k-1$ whole islands each, the remaining multicake contains 
at most $m'$ islands. Step \ref{step:complete} guarantees that it contains exactly $m'$ islands.
The remaining multicake has to be divided among $n' = n - \ell$ agents.
Since $n'\cdot k - n' + 1 = 
(n-\ell)(k-1)+1 = (nk-n+1)-\ell(k-1) \leq m - \ell(k-1) = m'$,
the integers $m',n',k$ still satisfy the input requirements of the algorithm.

For each of the remaining $n'$ agents, 
the allocated pieces are worth together at most $k\cdot \ell$, so the remaining multicake is worth for them at least:
\begin{align*}
(m+n-1) - k\ell
= (m-(k-1)\ell) + (n-\ell) - 1
= m' + n' - 1
\end{align*}
so the input requirements of algorithm \ref{alg:main} are satisfied.
Therefore, by the induction assumption, the recursive call solves problem \ref{prob:multicake-normalized} for the remaining agents.
\end{proof}

The following sections are devoted to solving the Partial Allocation problem. For presentation purposes, we first solve the problem for $k=1$ and $k=2$ and then present the solution for an arbitrary $k$, since each increase in $k$ requires a new technique.

\subsection{One Piece Per Agent}
\label{sub:k=1}
The following algorithm solves the partial allocation problem for $k=1$.
\begin{enumerate}
\item Let $B := \bestk{1}(C) = $ the island most valuable to agent \#1.
\item Ask each agent $j\in[n]$ to specify a real number $d_j\geq 0 $ as follows:
\begin{itemize}
\item If $V_j(B)\geq 1$ then $d_j$ is a number satisfying $V_j(B[d_j]) = 1$;
\item If $V_j(B)<1$ then $d_j = \infty$.
\end{itemize}
\item Choose an $i\in[n]$ such that $d_i$ is minimal (break ties arbitrarily). Return $L := \{i\}$ and the piece $Z_i := B[d_i]$.
\end{enumerate}
Let us show that this algorithm indeed solves problem \ref{prob:partial}.

The input conditions for that problem imply that the average value per island is at least $(m+n-1)/m\geq 1$. Hence, the island chosen in step 1 satisfies $V_1(B)\geq 1$.
Hence, in step 2, the first case holds for agent \#1, so $d_1$ is finite. 
Hence, in step 3, the minimal $d_i$ is finite.
The piece $Z_i$ contains $0 = k-1$ whole islands, plus a leftmost  sub-interval of a single island $B$. Its value for agent $i$ is exactly $1=k$. By the minimality of $d_i$, its value for the remaining agents is at most $1$ (for agents with $d_j=\infty$, its value is strictly less than $1$).
Hence, the output conditions for problem \ref{prob:partial} are satisfied.

\begin{remark}
In case there is $m=1$ island, the above algorithm reduces to the recursive step in the Last Diminisher algorithm \citep{Steinhaus1948Problem}.

The case $m\geq 1$ and $k=1$ was solved using a different algorithm at Example 4.2 in 
\citet{SegalHalevi2017Fair} and Lemma 11 in  \citet{segalhalevi2018redividing}.
A naive way to extend this algorithm to any $k>1$ is to replace each agent with $k$ clones having the same value function, run the algorithm for $k=1$, and then assign to each agent the union of the $k$ intervals allocated to his/her clones. This yields a value-guarantee of $k/(nk+m-1)$, which is much smaller than $k/(n+m-1)$ when $k$ is large.
To attain the high bar set by Theorem \ref{thm:prop}, we will need some new tools.
\end{remark}

\subsection{Two Pieces Per Agent}
\label{sub:k=2}
As a warm-up for the general case, this section develops a solution to the Partial Allocation Problem when $k=2$. 
At first glance, it seems that we could proceed as in the algorithm for $k=1$:
\begin{enumerate}
\item Let $B, B' := \bestk{1}(C) = $ the two best islands of agent \#1; note that $V_1(B\cup B')\geq 2$.
\item Ask each agent $j\in[n]$ to specify a $d_j\geq 0$ such that 
$V_j(B\cup B'[d_j]) = 2$;
\item Choose an $i$ such that $d_i$ is minimal, return $L = \{i\}$ and $Z_i := B\cup B'[d_i]$.
\end{enumerate}
But there is a problem: if there is an agent $j$ for whom $V_j(B)>2$, then this agent won't be able to find an appropriate $d_i$. Indeed, this agent will not agree that the whole island $B$ be given to any other agent, since its value for him is more than $k$.
Therefore, the above scheme can be followed only if there exists some island that \emph{all} $n$ agents value at less than $2$. We call such an island \emph{barren}. 
\begin{definition}
A subset $A\subseteq C$ is called \emph{barren}
if $V_j(A)< k$  for all agents $j\in[n]$.
\end{definition}
To handle this case we introduce Algorithm \ref{alg:mark}. It is termed a \emph{Mark Auction}:%
\footnote{
A mark auction 
is essentially a 
discrete version of the moving-knife procedure \citep{Dubins1961How}, 
or a simultaneous version of the last-diminisher procedure \citep{Steinhaus1948Problem}.
In last-diminisher, agent \#1 marks a piece of cake, each consecutive agent is allowed to mark a subset of the currently marked piece, and the last agent who marked a subset (i.e. the agent who marked the smallest piece) wins.
In a mark auction all marks are simultaneous,
and in addition, the marks are constrained to a pre-specified island. Constrained mark auctions are also useful for dividing a 2-dimensional cake  \citep{SegalHalevi2017Fair}.
}
Each agent ``bids'' by marking a piece of the cake, and the agent who has marked the smallest piece wins. 
The algorithm of subsection \ref{sub:k=1}
 is a special case of Algorithm \ref{alg:mark} with $k=1$ and $A=\emptyset$; the correctness proofs of both algorithms are similar too.

\begin{algorithm}[t]
\caption{
\label{alg:mark}
Mark Auction
}
\begin{algorithmic}[1]
\REQUIRE ~\\
\begin{itemize}
\item A \emph{barren} subset $A\subseteq C$ containing exactly $k-1$ whole islands.
\item A single island $B \subseteq C$ such that, for some agent $j\in[n]$, $V_j(A\cup B)\geq k$.
\end{itemize}

\ENSURE
A solution to the partial allocation problem \eqref{prob:partial}: an $i\in[n]$ and a piece $Z_i$ with:
\begin{itemize}
\item $V_i(Z_i)\geq k$;
\item $V_j(Z_i)\leq k$ for all $j\neq i$;
\end{itemize}
where $Z_i$ contains $k-1$ whole islands and a leftmost subset of a $k$-th island. 

\dottedline{}

\STATE
For each agent $j\in[n]$, define $d_j$ as follows:
\renewcommand{\theenumi}{\alph{enumi}}
\begin{enumerate}
\item If $V_j(A ~ \cup ~ B)<k$, then $d_j := \infty$;
\item Otherwise, $V_j(A ~ \cup ~ B)\geq k$: choose $d_j\geq 0$ such that $V_j(A ~ \cup ~ B[d_j]) = k$.
\end{enumerate}

\STATE
Choose an $i$ such that $d_i$ is minimal (break ties arbitrarily). 

\STATE
Return $L := \{i\}$ and the piece $Z_i := A~\cup~ B[d_i]$.
\end{algorithmic}
\end{algorithm}

\begin{lemma}
\label{lem:mark}
Whenever the input conditions of Algorithm \ref{alg:mark} are satisfied, it solves the partial allocation problem \ref{prob:partial} for any $k\geq 1$.
\end{lemma}
\begin{proof}
We first show that $d_j$ in step 1(b) is well-defined for all $j\in[n]$.
By the condition that $A$ is barren, $V_j(A)<k$. In step 1(b), $V_j(A) + V_j(B)\geq k$. Hence, 
there is a positive fraction of $B$ that, when it is added to $A$, the total value is exactly $k$. 

The input condition on $B$ implies that case 1(b) holds for at least one agent $j$. Hence, for this $j$, the $d_j$ is finite. 
Hence, the minimum $d_i$ found in step 2 is finite too.

The piece $Z_i$ contains, by construction, $ k-1$ whole islands ($A$), plus a leftmost sub-interval of length $d_i$ of a single island ($B$). Its value for agent $i$ is exactly $k$. By minimality of $d_i$, the value of $Z_i$ for the remaining agents is at most $k$ (for agents with $d_j=\infty$ its value is strictly less than $k$).
Hence, the output conditions for problem \ref{prob:partial} are satisfied.
\end{proof}

To use the Mark Auction with $k=2$, we need a barren island. What if there is no barren island? To handle this case we need a different tool. We use the concept of \emph{envy-free matching} \citep{Luria2013EnvyFree}.
Let $G = (X \cup Y, E)$ be a bipartite graph.
Recall that a \emph{matching} in $G$ is a subset $M\subseteq E$ such that each vertex in $X \cup Y$ is adjacent to at most one edge in $M$. We denote by $X_M$ ($Y_M$) the vertices of $X$ ($Y$) that are adjacent to an edge of $M$.
\begin{definition}
\label{def:efm}
Let $G = (X \cup Y, E)$ be a bipartite graph.
A matching $M\subseteq E$ is called 
an \emph{envy-free matching} (EFM)%
\footnote{
Envy-free matching is unrelated to 
\emph{envy-freeness for mixed goods}, 
which is also abbreviated as EFM \citep{bei2019fair}.
}
if no unmatched vertex in $X$ is adjacent to a matched vertex in $Y$. Formally:
\begin{align*}
\forall x\in X\setminus X_M: ~~ \forall y\in Y_M: ~~ (x,y)\notin E.
\end{align*}
\end{definition}
In an EFM, an unmatched $x$ does not ``envy'' any matched $x'$, because it does not ``like'' any matched $y'$ anyway.%
\footnote{
The above definition of envy-free matching involves an unweighted graph, where each ``agent'' (a vertex in $X$) either likes or dislikes each ``object'' (a vertex in $Y$). 
Some other papers, such as \citet{gan2019envy}, study an envy-free assignment in a weighted graph, in which each agent has a ranking on the objects.
}
Any $X$-perfect matching (i.e., a matching in which every vertex in $X$ is matched) is envy-free, and the empty matching is envy-free too. The following theorem provides a sufficient condition for the existence of a nonempty EFM (here $N_G(X)$ is the set of neighbors of $X$ in the graph $G$):
\begin{efmtheorem*}
Let $G = (X\cup Y, E)$ be a bipartite graph. If $|N_G(X)|\geq |X|\geq 1$, then $G$ admits a nonempty envy-free matching, and such a matching can be found by a polynomial-time algorithm.
\end{efmtheorem*}
The existential part of the EFM theorem was proved by \citet{Luria2013EnvyFree}.
Recently,   \citet{aigner2019envy} presented a polynomial-time algorithm for finding an EFM of maximum cardinality; in particular, it finds a nonempty EFM whenever it exists.
To make this paper self-contained, a short proof of the EFM theorem is presented in Appendix \ref{sec:efm}.%
\footnote{
\citet{amanatidis2017approximation} also presented a polynomial-time algorithm for finding a nonempty EFM when the conditions of the EFM theorem are satisfied. Their algorithm was stated not for general graphs but for a specific graph constructed in the context of fair allocation of indivisible goods. 
}

We now present an algorithm for Partial Allocation with $k=2$.
Consider two main cases.

\noindent
\textbf{Case \#1:} There exists a barren island $A_0$.

For each agent $j\in[n]$, let $B_j$ be the island most valuable to $j$ (breaking ties arbitrarily). Consider two sub-cases:

\noindent
\emph{Case \#1.1:} $V_j(A_0\cup B_j)\geq 2$ for at least one agent $j$. Then:
\begin{enumerate}
\item Do a Mark Auction with $A=A_0$ and $B=B_j$.
\end{enumerate}

\noindent
\emph{Case \#1.2:} For all agents $j$, $V_j(A_0\cup B_j)< 2$. Then:
\begin{enumerate}
\item Let $B_1, B_1'$ be the two islands most valuable for agent \#1.
\item Do a Mark Auction with $A=B_1$ and $B=B_1'$.
\end{enumerate}

\noindent
\textbf{Case \#2:} There is no barren island.
\begin{enumerate}
\item Construct a bipartite graph $G = (X\cup Y,E)$ where $X$ is the set of agents, $Y$ the set of islands, and there is an edge between agent $i$ and island $B$ whenever $V_i(B)\geq 2$.
\item Find a nonempty envy-free matching in $G$. 
\item Return $L = $ the matched agents and $(Z_i)_{i\in L}$  = the matched pieces.
\end{enumerate}

\begin{claim}
The algorithm described above solves the partial allocation problem \ref{prob:partial} when $k=2$.
\end{claim}
\begin{proof}
We consider each case in turn.

\emph{Case \#1.1:} Since $A_0$ is barren, and $V_j(A_0\cup B_j)\geq k$ for at least one $j\in[n]$, the input requirements of the Mark Auction are satisfied, so by Lemma \ref{lem:mark} it finds a partial allocation.

\emph{Case \#1.2:}
For all agents $j\in[n]$,
we have $V_j(A_0\cup B_j)<2$, where $B_j$ is the most valuable island of $j$. Hence, 
every agent values every island at less than $2$ --- all islands are barren.
In particular, $B_1$ is barren.
On the other hand, $V_1(B_1\cup B_1')\geq 2$,
since the average value per island is at least $1$.
Again the input requirements of the Mark Auction are satisfied, so by Lemma \ref{lem:mark} it finds a partial allocation.

\emph{Case \#2:} 
Every island is valued at least $2$ by at least one agent. Hence, in the bipartite graph $G$ constructed in step 1, each island in $Y$ has at least one neighbor in $X$. Therefore, $|N_G(X)| = m\geq n+1 \geq |X|$. The EFM theorem implies that $G$ admits a nonempty envy-free matching.
Each matched piece contains $1 = k-1$ whole island, and its value for its matched agent is at least $2=k$.
Because the matching is envy-free, each matched piece is not adjacent to each unmatched agent, so the value of each matched piece to an unmatched agent is less than $2$. Hence the output requirements of problem \ref{prob:partial} are satisfied and a partial allocation is found.
\end{proof}

\subsection{Many Pieces per Agent}
\label{sub:k=k}
This section presents a partial-allocation algorithm for any $k\geq 1$.
The general scheme is similar to the case $k=2$: we find either an appropriate input to a mark auction, or an appropriate input to an envy-free matching. 
We call a pair satisfying the requirements of a mark auction, a \emph{threshold pair}.
\begin{definition}
A \emph{threshold pair} is a pair $(A,B)$ such that:
\begin{itemize}
\item $A\subseteq C$ is a barren subset of $k-1$ islands;
\item $B\subseteq C$ is a single island such that, for some agent $j\in[n]$, $V_j(A\cup B)\geq k$.
\end{itemize}
\end{definition}
Given a threshold pair, the Mark Auction (algorithm \ref{alg:mark}) can be used to find a partial allocation.
Below it is shown that, given a single barren subset of $k-1$ islands, a threshold pair can be found.

\begin{algorithm}[t]
\caption{
\label{alg:threshold}
Finding a threshold pair.
}
\begin{algorithmic}[1]
\REQUIRE A barren subset $A_0\subseteq C$ containing exactly $k-1$ whole islands.

\ENSURE A threshold pair $(A^*, B^*)$, that is:
\begin{itemize}
\item $A^*\subseteq C$ is a barren subset of $k-1$ islands;
\item $B^*\subseteq C$ is a single island such that, for some agent $j\in[n]$, $V_j(A^*\cup B^*)\geq k$.
\end{itemize}

\dottedline{}

\STATE For every agent $j\in[n]$, let $B_j$ be
a set of $k-|A_0|$ islands in $C\setminus A_0$ most valuable to $j$.
\STATE If $V_j(A_0\cup B_j)\geq k$ for some $j\in[n]$, 
then terminate and return:
\begin{itemize}
\item $B^* :=$ an arbitrary whole island from $B_j$.
\item $A^* := A_0 \cup B_j \setminus B^*$.
\end{itemize}

\STATE Else, $V_j(A_0\cup B_j) < k$ for all $j\in[n]$:
\begin{itemize}
\item Remove an arbitrary island from $A_0$.
\item Go back to step 1.
\end{itemize}
\end{algorithmic}
\end{algorithm}

\begin{lemma}
\label{lem:threshold}
Given any barren subset of $k-1$ islands, Algorithm \ref{alg:threshold} finds a threshold pair.
\end{lemma}

Before proving Lemma \ref{lem:threshold} formally, let us describe informally how Algorithm \ref{alg:threshold} works.

We are given $A_0$ --- a barren subset of $k-1$ islands.
We first look for an island $B_j\subseteq C\setminus A_0$ such that  $V_j(B_j)\geq k-V_j(A_0)$ for some agent $j\in[n]$.
If such a $B_j$ exists, 
then $(A_0,B_j)$ is a threshold pair, since $V_j(A_0\cup B_j)\geq k$, so we are done.

If such a $B_j$ does not exist, we know that \emph{every} subset of $k$ islands that contains the $k-1$ islands of $A_0$ is barren. 
Every barren subset obviously remains barren when any island is removed from it. 
Specifically, if we remove an arbitrary island $A'$ from $A_0$, then every subset of $k-1$ islands that contains the $k-2$ islands remaining in $A_0$ is barren: it is barren even with $A'$, so it is certainly barren without $A'$.
Now, we look for a pair of islands $B_j,B_j'$ such that $V_j(B_j'\cup B_j)\geq k-V_j(A_0)$ for some agent $j\in[n]$.
If such a pair exists, then 
$(A_0\cup B_j', ~B_j)$ is a threshold pair,
since, by the previous paragraph, $A_0\cup B_j'$ is a barren set of $k-1$ islands, and 
$V_j(A_0\cup B_j'\cup B_j)\geq k$.

If such a pair $B_j,B_j'$ does not exist, then we know that \emph{every} subset of $k$ islands that contains the $k-2$ islands of $A_0$ is barren.
We remove another island from $A_0$; now every subset of $k-1$ islands that contains the $k-3$ islands remaining in $A_0$ is barren. We look for a triplet $B_j,B_j',B_j''$ such that $V_j(B_j'' \cup B_j'\cup B_j)\geq k-V_j(A_0)$ for some agent $j$.
If such a triplet exists, then $(A_0\cup B_j''\cup B_j',~B_j)$ is a threshold pair, since $A_0\cup B_j''\cup B_j'$ is a barren set of $k-1$ islands.

If such a triplet does not exist, then we keep removing islands from $A_0$, one by one, and looking for some $k-|A_0|$ islands that, for some agent $j$, are worth at least $k-V_j(A_0)$.
This process must end at some point, since the average value per island is at least $1$, so when $A_0$ becomes empty, for \emph{every} agent, the $k$ most valuable islands are worth at least $k$.

We now formally prove the correctness of Algorithm \ref{alg:threshold}.

\begin{proof}[Proof of Lemma \ref{lem:threshold}]
We first show that the algorithm maintains the following invariant:
\begin{align}
\label{eq:conserve}
\tag{@}
\text{\emph{At step 1, every set of $k-1$ islands that contains $A_0$ is barren.}}
\end{align}

This holds at the first time step 1 is reached, since by the input condition, $A_0$ contains $k-1$ islands and it is barren.
If the algorithm does not terminate at step 2, it means that every set of $k$ islands that contains $A_0$ is barren. 
At step 3, an arbitrary island (say, $A'$) is removed from $A_0$.
Now, every set of $k-1$ islands that contains the islands remaining in $A_0$ is barren: it is barren even with $A'$, so it is certainly barren without $A'$. Hence \eqref{eq:conserve} remains valid.

Next, we show that the algorithm must terminate. 
Indeed, at each iteration, $A_0$ shrinks by one island, so 
at the $k$-th iteration, $A_0$ is empty. Then, for every $j\in[n]$, $B_j$ is a set of $k$ islands most valuable to $j$. Since the average value per island is at least $(m+n-1)/m\geq 1$, $V_j(B_j)\geq k$ for every $j\in[n]$, so the algorithm terminates at step 2.

Suppose the algorithm terminates at step 2 of iteration $t$, for some $t\in\{1,\ldots,k\}$.
Then, $A_0$ contains $k-t$ islands and $B_j$ contains $t$ islands. Therefore, $A^*$ contains $(k-t)+(t)-1 = k-1$ islands.
$A^*$ contains $A_0$, so by the invariant \eqref{eq:conserve}, $A^*$ is barren. 
By the termination condition, $V_j(A^*\cup B^*) = V_j(A_0\cup B_j)\geq k$. Hence, $(A^*,B^*)$ is a threshold pair as claimed.
\end{proof}

\begin{algorithm}[t]
\caption{
\label{alg:partial}
Partial Allocation Algorithm for any $k\geq 1$.
}
\begin{algorithmic}[1]
\REQUIRE ~\\
\begin{itemize}[noitemsep]
\item Positive integers $m, n, k$ with $m\geq nk-n+1$;
\item A multicake $C$ with $m$ islands; 
\item $n$ value-measures $V_1,\ldots,V_n$ on $C$ having $V_i(C) \geq n+m-1$ for all $i\in[n]$.
\end{itemize}

\ENSURE A nonempty set of agents $L\subseteq [n]$, and a set of pieces $(Z_i)_{i\in L}$ with for all $i\in L$:
\begin{itemize}
\item The piece $Z_i$ contains exactly $k-1$ whole islands of $C$, plus possibly the leftmost part of some $k$-th island.
\item $V_i(Z_i)\geq k$.
\item $V_j(Z_i)\leq k$ for all $j\not\in L$.
\end{itemize}
\dottedline{}

\STATE Arbitrarily select $n$ pairwise-disjoint subsets of $k-1$ islands each.

\STATE
\textbf{Case \#1:} At least one of these $n$ subsets is barren. 
\begin{itemize}
\item Using Algorithm \ref{alg:threshold}, find a threshold pair $A^*, B^*$.
\item Using Algorithm \ref{alg:mark} (mark auction) with $A^*, B^*$, find a partial allocation.
\end{itemize}

\STATE
\textbf{Case \#2:} None of these $n$ subsets is barren. 
\begin{itemize}
\item Construct a bipartite graph $G = (X\cup Y,E)$ where $X$ is the set of $n$ agents, $Y$ the set of $n$ subsets of $k-1$ islands, and there is an edge between agent $i$ and subset $B_j$ whenever $V_i(B_j)\geq k$.
\item Find a nonempty envy-free matching in $G$. 
\item Return $L = $ the matched agents and $(Z_i)_{i\in L}$  = the matched pieces.
\end{itemize}

\end{algorithmic}
\end{algorithm}

Algorithm \ref{alg:partial} uses Algorithm \ref{alg:threshold} 
to solve the partial allocation problem for any $k\geq 1$.
\begin{lemma}
\label{lem:partial}
Algorithm \ref{alg:partial} solves the Partial Allocation problem \ref{prob:partial}.
\end{lemma}
\begin{proof}
By the input conditions, the number of islands is $m\geq nk-n+1 > n\cdot (k-1)$. Hence there are indeed at least $n$ pairwise-disjoint subsets of $k-1$ islands, as required by the first step.

In Case \#1, Lemma \ref{lem:threshold} implies that Algorithm \ref{alg:threshold} indeed finds a threshold pair, and 
Lemma \ref{lem:mark} implies that Algorithm \ref{alg:mark} indeed finds a partial allocation.

In Case \#2, since each of the $n$ subsets is not barren, each vertex of $Y$ is adjacent to at least one vertex of $X$. Hence $|N_G(X)| = |Y| = |X| = n \geq 1$. By the EFM theorem, a nonempty envy-free matching exists. 
By the envy-freeness of the matching, the resulting $L$ and $(Z_i)_{i\in L}$ are indeed a partial allocation.
\end{proof}

Now we can finally prove the main theorem.
\begin{proof}[Proof of Theorem \ref{thm:prop}]
Lemma \ref{lem:upperbound} proves the upper bound $r(m,k,n)\leq \min\left({1\over n}, {k\over m+n-1}\right)$.

Lemma \ref{lem:partial} implies that Algorithm \ref{alg:partial} solves the Partial Allocation Problem.
By plugging this algorithm into Algorithm \ref{alg:main}, 
Lemma \ref{lem:partial-implies-multicake} now implies that 
Algorithm \ref{alg:main} solves the Normalized Multicake Division Problem. By Lemma \ref{lem:normalization}
it solves the Multicake Division Problem. Hence 
$r(m,k,n)\geq \min\left({1\over n}, {k\over m+n-1}\right)$.
\end{proof}

\begin{remark}
\label{rem:cuts}
Algorithm \ref{alg:partial} makes at most one cut on a single island: in case \#1 it makes one cut during the mark auction, and in case \#2 it makes no cuts at all.
Therefore, Algorithm \ref{alg:main}  makes at most $n-1$ cuts on the multicake.
In some cases $n-1$ cuts may be necessary, e.g. when $m=1$.
\end{remark}

\subsection{Example}
\label{sub:example}
The operation of the multicake division algorithm is illustrated below.

In this example, $k=3$ and there are $n=4$ agents. Initially there are $m=7$ islands. 
Since $m < nk-n+1$, two dummy islands are added as in Lemma \ref{lem:normalization}, and the total number of islands becomes $9$. 
The agents' valuations are normalized such that the value of the entire multicake is $9+4-1=12$. The normalized values are shown below, where $C_j$ denotes the $j$-th island (in an arbitrary order).

\begin{center}
\begin{filecontents*}{example11.csv}
agent;ca;cb;cc;cd;ce;cf;cg;ch;ci;total
1;5;2;1;1;1;1;1;0;0;12
2;2;4;1;1;1;1;2;0;0;12
3;1;1;1;2;2;2;3;0;0;12
4;1;1;1;1;1;1;6;0;0;12
\end{filecontents*}
\csvreader
[tabular=||c||c|c|c|c|c|c|c|c|c||c||,
table head=\hline Agent & $C_1$ &  $C_2$ &  $C_3$ &  $C_4$ &  $C_5$  &  $C_6$ &  $C_7$  &  $C_8$ &  $C_9$ & Sum \\\hline,
late after line=\\, 
late after last line=\\\hline,
/csv/separator=semicolon,
head to column names]{example11.csv}{}%
{
\agent & \ca &
\cb & \cc &
\cd & \ce &
\cf & \cg &
\ch & \ci & \total
}
\end{center}

Algorithm \ref{alg:partial} is used to find a partial allocation.
The $n$ subsets of $2$ islands are arbitrarily chosen as: $C_1\cup C_2$, $C_3\cup C_4$, $C_5\cup C_6$, $C_7\cup C_8$. 
Each of these is worth at least $3$ to some agent --- no pair is barren. Therefore we look for an envy-free matching. The bipartite graph is:

\begin {center}
\begin {tikzpicture}[-latex ,auto ,node distance =2 cm and 3cm ,on grid, semithick , state/.style ={ circle ,top color =white , bottom color=blue!20, draw,blue , text=blue , minimum width =1 cm}]
a\node[state] (A1) {$1$};
\node[state] (A2) [right=of A1] {$2$};
\node[state] (A3) [right=of A2] {$3$};
\node[state] (A4) [right=of A3] {$4$};
\node[state] (Y1) [below=of A1] {$C_1 \cup C_2$};
\node[state] (Y2) [below=of A2] {$C_3 \cup C_4$};
\node[state] (Y3) [below=of A3] {$C_5 \cup C_6$};
\node[state] (Y4) [below=of A4] {$C_7 \cup C_8$};
\path (A1) edge (Y1);
\path (A2) edge (Y1);
\path (A3) edge (Y2);
\path (A3) edge (Y3);
\path (A3) edge (Y4);
\path (A4) edge (Y4);
\end{tikzpicture}
\end{center}
It contains a nonempty envy-free matching, for example matching 
agent $4$ to $C_7\cup C_8$
and agent $3$ to $C_5\cup C_6$. Agents $3$ and $4$ now receive two whole islands worth for them at least $3$ and are removed from further consideration. 
The remaining valuations are:

\begin{center}
\begin{filecontents*}{example12.csv}
agent;ca;cb;cc;cd;ce;cf;cg;ch;ci;total
1;5;2;1;1;1;1;1;0;0;9
2;2;4;1;1;1;1;2;0;0;8
\end{filecontents*}
\csvreader
[tabular=||c||c|c|c|c|c||c||,
table head=\hline Agent & $C_1$ &  $C_2$ &  $C_3$ &  $C_4$ &    $C_9$ & Sum \\\hline,
late after line=\\, 
late after last line=\\\hline,
/csv/separator=semicolon,
head to column names]{example12.csv}{}%
{
\agent & \ca &
\cb & \cc &
\cd & \ci & \total
}
\end{center}

Again Algorithm \ref{alg:partial} is used to find a partial allocation.
The $n$ subsets of $2$ islands are arbitrarily chosen as: $C_1\cup C_2$ and $C_3\cup C_4$. 
The second pair is barren --- it is worth less than $3$ to both agents. Therefore,
Algorithm \ref{alg:threshold} is used to find a threshold pair.
Initially $A_0 = C_3\cup C_4$. 
$B_1$ is set to $C_1$ --- the most valuable island for agent 1.
Similarly, $B_2 = C_2$.
Since $V_1(A_0\cup B_1)\geq k$, the algorithm returns $A^* := C_3\cup C_4$ and $B^* := C_1$.
Now a Mark Auction (algorithm \ref{alg:mark}) is used: each agent $j\in\{1,2\}$ is asked to mark some $d_j\geq 0$ such that $V_j(C_3\cup C_4\cup C_1[d_j])=k$. For the sake of example, suppose that agent 1 marks $d_1=0.2$ and agent 2 marks $d_2=0.5$.
Then agent 1 wins and receives two whole islands ($C_3$ and $C_4$), plus an interval of length $0.2$ at the left of $C_1$. The remaining values are (where $C_1'$ is the remaining part at the right of $C_1$):

\begin{center}
\begin{filecontents*}{example13.csv}
agent;ca;cb;cc;cd;ce;cf;cg;ch;ci;total
2;1.6;4;1;1;1;1;2;0;0;5.6
\end{filecontents*}
\csvreader
[tabular=||c||c|c|c||c||,
table head=\hline Agent & $C_1'$ &  $C_2$ &   $C_9$ & Sum \\\hline,
late after line=\\, 
late after last line=\\\hline,
/csv/separator=semicolon,
head to column names]{example13.csv}{}%
{
\agent & \ca &
\cb & \ci & \total
}
\end{center}

By the termination condition of Algorithm \ref{alg:main}, agent 2 receives the three remaining islands.

\subsection{Open questions}
\label{sub:algorithm-open}
For a fixed $k$, our Multicake Division Algorithm attains the highest possible value-guarantee. This raises the question of what happens if different agents have different values of $k$, corresponding to different preferences about the number of pieces.
Suppose some agents insist on getting a small number of pieces even if it entails a smaller value-guarantee, while others insist on getting a large value even if it requires many pieces. Is it possible to give a personalized guarantee to each agent? 
\begin{open}
Given positive integers $k_i$ for all $i\in[n]$,
is it possible to allocate to each agent $i$ a piece made of at most $k_i$ intervals,
with a value of at least  $\min\left({1\over n}, {k_i\over m+n-1}\right)\cdot V_i(C)$?
\end{open}

Instead of the the hard constraint of having at most $k$ pieces per agent,
one can consider a softer model where the number of pieces decreases the value by some known amount. 
For example, suppose that managing pieces on $k$ different islands incurs an expense of $E(k)$. Then the net utility of agent $i$ from an allotment $Z_i$ overlapping $k_i$ islands is $V_i(Z_i)-E(k_i)$. 
\begin{open}
What value-guarantees are attainable in that model?
\end{open}

\section{Relative Value Guarantee}
\label{sec:relprop}
This section aims to prove Theorem \ref{thm:relprop} by giving each agent at least $1/n$ of the value of his/her best $k$ islands.
See Problem \ref{prob:rel-multicake} below.
As in the previous section, it is convenient to work with the normalized format specified by Problem \ref{prob:rel-multicake-normalized} below.

\begin{problem}
\caption{
\label{prob:rel-multicake}
Multicake Division Problem --- Relative Value Guarantee
}
\begin{algorithmic}
\REQUIRE ~\\
\begin{itemize}
\item Positive integers $m, n, k$; 
\item A multicake $C$ with $m$ islands;  
\item $n$ value-measures $V_1,\ldots,V_n$ on $C$.
\end{itemize}

\ENSURE A $k$-feasible allocation $\mathbf{Z} = Z_1,\ldots,Z_n$ in which for all $i\in[n]$:
\begin{align*}
{V_i(Z_i)} \geq {1\over n} \cdot  \vbestk{i}{C}.
\end{align*}
\end{algorithmic}
\end{problem}

\begin{problem}
\caption{
\label{prob:rel-multicake-normalized}
Normalized Multicake Division Problem --- Relative Value Guarantee 
}
\begin{algorithmic}
\REQUIRE ~\\
\begin{itemize}[noitemsep]
\item Positive integers $m, n, k$ with $m\geq nk-n+1$;
\item A multicake $C$ with $m$ islands; 
\item $n$ value-measures $V_1,\ldots,V_n$ on $C$, such that for every agent $i\in[n]$:
\begin{itemize}
\item At most $k$ islands have a positive value; 
\item At least $m-k$ islands have a zero value;
\item $V_i(C) \geq k n$.
\end{itemize}
\end{itemize}

\ENSURE A $k$-feasible allocation $\mathbf{Z} = Z_1,\ldots,Z_n$ in which for all $i\in[n]$:
\begin{align*}
V_i(Z_i)\geq k.
\end{align*}
\end{algorithmic}
\end{problem}

\begin{lemma}
\label{lem:rel-normalization}
Any instance of problem \ref{prob:rel-multicake} can be reduced to an instance of problem \ref{prob:rel-multicake-normalized}.
\end{lemma}
\begin{proof}
Given an instance of problem \ref{prob:rel-multicake}, construct an instance of problem \ref{prob:rel-multicake-normalized} 
with the same $n,k$ and with $m' = \max(m, nk-n+1)$;
if $m<nk-n+1$, then add dummy islands whose value for all agents is $0$, such that the total number of islands becomes $nk-n+1$.
Construct value-measures $V_1',\ldots,V_n'$ as the integrals of the following value-density functions:
\begin{align*}
v_i'(x) &:= 
\begin{cases}
\dfrac{k \cdot n}{\vbestk{i}{C}}
v_i(x)
&
x \in \bestk{i}(C)
\\
0
&
\text{otherwise}
\end{cases}
\\
V_i'(Z_i) &:= \int_{x\in Z_i} v_i'(x) dx.
\end{align*}
For all $i\in[n]$, $V_i'(C) = \dfrac{k \cdot n}{\vbestk{i}{C}}\cdot V_i\left(\bestk{i}(C)\right) = k\cdot n$. Moreover, only the (at most) $k$ best islands have a positive value by $V_i'$. Hence, these valuations satisfy the input requirements of problem \ref{prob:rel-multicake-normalized}.

The solution is a $k$-feasible allocation $\mathbf{Z}$ such that $V_i'(Z_i)\geq k$ for all $i\in[n]$. 
By normalization, 
$
{V_i(Z_i)} \geq k\cdot \dfrac{\vbestk{i}{C}}{k\cdot n} = {1\over n} \cdot  \vbestk{i}{C}$,
so the same allocation $\mathbf{Z}$ solves problem \ref{prob:rel-multicake}.
\end{proof}

Problem \ref{prob:rel-multicake-normalized} can be solved by the same algorithm used in the previous section (Algorithm \ref{alg:main}) --- only the normalization should be different.
To prove this claim we follow the proof steps of Section \ref{sec:algorithm}.

\begin{lemma}
Given a correct algorithm for solving the Partial Allocation problem \eqref{prob:partial}, 
when Algorithm \ref{alg:main} is executed with valuations satisfying the input conditions of Problem \ref{prob:rel-multicake-normalized}, 
it correctly solves Problem \ref{prob:rel-multicake-normalized}.
\end{lemma}
\begin{proof}
By induction on $n$. If $n=1$ then the algorithm terminates at step 1. The agent values positively at most $k$ islands, and values them at least $k\cdot 1 = k$. Therefore the $k$ most valuable islands are worth for this agent at least $k$.

Suppose now that $n>1$
and that the claim holds whenever there are less than $n$ agents.
At step 2, at least one agent receives a piece, and the partial allocation is $k$-feasible.
Each allocated piece is worth at least $k$ to its receiver.

Let $\ell = |L|$ and $n' := n-\ell$.
For the remaining $n'$ agents,
the allocated pieces are worth together at most $k\cdot \ell$,
so the remaining multicake is worth at least $k\cdot n'$. 
Since the partial allocation does not create new islands with a positive value, at most $k$ islands have a positive value for each remaining agent. 
Hence, the input requirements of problem \ref{prob:rel-multicake-normalized} are satisfied.
By the induction assumption, the recursive call solves problem \ref{prob:rel-multicake-normalized} for the remaining agents.
\end{proof}

Algorithm \ref{alg:mark} --- the Mark Auction --- does not depend in any way on the normalization, so it still finds a partial allocation whenever its input conditions are satisfied (Lemma \ref{lem:mark}).

Algorithm \ref{alg:threshold} still finds a threshold pair; the proof is almost identical to Lemma \ref{lem:threshold}.
The only part of the proof that depends on the normalization
is the claim that the algorithm must terminate when $A_0$ is empty (the before-last paragraph of the proof).
To prove that the algorithm terminates, 
we have to prove that, for every agent, the value of the $k$ best islands is at least $k$.
With the normalization of Problem \ref{prob:multicake-normalized}, 
this follows from the fact that the total multicake value is at least $m$, so the average value per island is at least $1$.
With the normalization of Problem \ref{prob:rel-multicake-normalized},
the total multicake value is $nk$ which may be smaller than $m$, so the average value per island may be less than $1$.
However, the condition that only the $k$ best islands have a positive value guarantees that the value of these islands is at least $n k$, which is at least $k$.

Similarly, Algorithm \ref{alg:partial} still finds a partial allocation
when the input conditions are replaced with those of Problem \ref{prob:rel-multicake-normalized}; the proof is identical to Lemma \ref{lem:partial}.

Combining all the above arguments proves Theorem \ref{thm:relprop}.

%

\begin{remark}
All the proofs regarding the value guarantees use only the individual value functions. Therefore, when using Algorithm \ref{alg:main}, 
it is possible to let different agents use different normalizations simultaneously and attain their more favorable guarantee. In particular, an agent for whom the absolute guarantee of Theorem \ref{thm:prop} is higher may use the normalization of Problem \ref{prob:multicake-normalized}, while an agent for whom the absolute guarantee of Theorem \ref{thm:relprop} is higher may use the normalization of Problem \ref{prob:rel-multicake-normalized}.
An example is given below.
\end{remark}

\subsection{Example}
\label{sub:example-rel}
In this example, $k=2$ and there are $n=4$ agents and $m=7$ islands. 
Initially, the agents' valuations are normalized such that the value of the entire multicake is $7+4-1=10$. The normalized values are shown below.
\begin{center}
\begin{filecontents*}{example21.csv}
agent;ca;cb;cc;cd;ce;cf;cg;total
1;6;3.6;0.4;0;0;0;0;10
2;2;4;1;1.6;0.4;0;1;10
3;0;2;3;1;2;1;1;10
4;1;1;1;0;1;4;2;10
\end{filecontents*}
\csvreader
[tabular=||c||c|c|c|c|c|c|c||c||,
table head=\hline Agent & $C_1$ &  $C_2$ &  $C_3$ &  $C_4$ &  $C_5$  &  $C_6$ &  $C_7$  &  Sum \\\hline,
late after line=\\, 
late after last line=\\\hline,
/csv/separator=semicolon,
head to column names]{example21.csv}{}%
{
\agent & \ca &
\cb & \cc &
\cd & \ce &
\cf & \cg &
\total
}
\end{center}
Theorem \ref{thm:prop} guarantees to each agent a value of at least $2$. 
However, for agent 1, Theorem \ref{thm:relprop} guarantees a higher value --- at least $9.6/4 = 2.4$. 
To implement this higher guarantee, the valuations of agent 1 are re-normalized as in Problem \ref{prob:rel-multicake-normalized}
(see Section \ref{sec:relprop}):

\begin{center}
\begin{filecontents*}{example22.csv}
agent;ca;cb;cc;cd;ce;cf;cg;total
1;5;3;0;0;0;0;0;8
2;2;4;1;1.6;0.4;0;1;10
3;0;2;3;1;2;1;1;10
4;1;1;1;0;1;4;2;10
\end{filecontents*}
\csvreader
[tabular=||c||c|c|c|c|c|c|c||c||,
table head=\hline Agent & $C_1$ &  $C_2$ &  $C_3$ &  $C_4$ &  $C_5$  &  $C_6$ &  $C_7$  &  Sum \\\hline,
late after line=\\, 
late after last line=\\\hline,
/csv/separator=semicolon,
head to column names]{example22.csv}{}%
{
\agent & \ca &
\cb & \cc &
\cd & \ce &
\cf & \cg &
\total
}
\end{center}
Now, all agents are given a value of at least $2$ as follows.

Algorithm \ref{alg:partial} is used to find a partial allocation.
The $n$ subsets of $1$ island are arbitrarily chosen as: $C_1$, $C_2$, $C_3$, $C_4$. 
The island $C_4$ is barren --- it is worth less than $2$ to all agents. Therefore,
Algorithm \ref{alg:threshold} is used to find a threshold pair, starting at $A_0 := C_4$. 
$B_1$ is set to $C_1$ --- the most valuable island for agent 1.
Since $V_1(A_0\cup B_1)\geq k$, the algorithm returns $A^* := C_4$ and $B^* := C_1$.
Now a Mark Auction is used: each agent $j$ is asked to mark some $d_j\geq 0$ such that $V_j(C_4\cup C_1[d_j])=2$.
Since for agents $3$ and $4$ the value of $C_4\cup C_1$ is less than $2$, they report $d_3 = d_4 = \infty$.
For the sake of example, suppose that the other marks are $d_1=0.4$ and $d_2=0.2$.
Then agent 2 wins and receives $C_4$  plus an interval of length $0.2$ at the left of $C_1$. The remaining values are (where $C_1'$ is the remainder at right of $C_1$):

\begin{center}
\begin{filecontents*}{example23.csv}
agent;ca;cb;cc;cd;ce;cf;cg;total
1;4;3;0;0;0;0;0;7
3;0;2;3;1;2;1;1;9
4;0.8;1;1;0;1;4;2;9.8
\end{filecontents*}
\csvreader
[tabular=||c||c|c|c|c|c|c||c||,
table head=\hline Agent & $C_1'$ &  $C_2$ &  $C_3$ &   $C_5$  &  $C_6$ &  $C_7$  &  Sum \\\hline,
late after line=\\, 
late after last line=\\\hline,
/csv/separator=semicolon,
head to column names]{example23.csv}{}%
{
\agent & \ca &
\cb & \cc &
 \ce &
\cf & \cg &
\total
}
\end{center}
Algorithm \ref{alg:partial} is used again, and let's assume that the $3$ subsets of $1$ island are $C_1'$ and $C_2$ and $C_3$. None of these islands is barren, so we find an envy-free matching in the following graph:

\begin {center}
\begin {tikzpicture}[-latex ,auto ,node distance =1.5 cm and 2cm ,on grid, semithick , state/.style ={ circle ,top color =white , bottom color=blue!20, draw,blue , text=blue , minimum width =1 cm}]
a\node[state] (A1) {$1$};
\node[state] (A3) [right=of A3] {$3$};
\node[state] (A4) [right=of A3] {$4$};
\node[state] (Y1) [below=of A1] {$C_1'$};
\node[state] (Y2) [below=of A3] {$C_2$};
\node[state] (Y3) [below=of A4] {$C_3$};
\path (A1) edge (Y1);
\path (A1) edge (Y2);
\path (A3) edge (Y2);
\path (A3) edge (Y3);
\end{tikzpicture}
\end{center}
The graph admits a nonempty envy-free matching, for example, matching agent 1 with $C_1'$ and agent $3$ with $C_3$.
Each of these two agents receives one island. 
The remaining agent $4$ receives his two most valuable remaining islands, namely $C_6$ and $C_7$.

\ifdefined\mms
\subsection{Open question}
In the study of fair allocation of indivisible items, 
a prominent fairness criterion is the \emph{maximin share}  \citep{budish2011combinatorial}.
Given an arbitrary set $C$ and $n\geq 1$ agents, the maximin share (MMS) of agent $i\in[n]$ is the largest value that $i$ can get by
constructing a feasible partition of $C$ into $n$ subsets and getting the least valuable part in that partition.
When $C$ is a multicake, 
``feasible'' means that each part contains at most $k$ intervals. 
So the MMS of agent $i$, denoted $\mmsk{i}(C)$, is the largest value that $i$ can get by
constructing a partition of (a subset of) $C$ into $n$ parts each of which contains at most $k$ intervals, and getting the least valuable part.
An allocation is \emph{MMS-fair} if the value allocated to each agent $i$ is at least $\mmsk{i}(C)$.

In the classing setting, in which $k = m = 1$, the MMS of all agents equals $1/n$ of the total cake value, so MMS-fairness is equivalent to proportionality. 
MMS-fairness becomes interesting when proportionality cannot be guaranteed, such as for a multicake with $m>k$.


Since each agent can ``clone'' itself $n$ times, Theorem \ref{thm:prop} implies that for all $i\in[n]$,
$
\mmsk{i}(C) \geq \min\left({1\over n},{k\over m+n-1}\right)\cdot V_i(C)
$.
Similarly, Theorem \ref{thm:relprop} implies that
$
\mmsk{i}(C) \geq {1\over n}\cdot V_i(\bestk{i}(C))
$.
But the MMS may be strictly larger than both these expressions. As an example, suppose there are $n=10$ agents and $m=10$ islands and $k=1$.
Suppose Alice values all islands uniformly at $1/10$ of the total value.
Then her MMS is $1/10$, 
but Theorem \ref{thm:prop} guarantees her only $1/19$ of the total value, while Theorem \ref{thm:relprop} guarantees her only
$1/100$.
This raises the question of whether an MMS-fair allocation always exists.

It is easy to show that Algorithm \ref{alg:main} always finds an MMS-fair allocation when $k=1$, but may fail to find an MMS-fair allocation when $k\geq 2$
(see Appendix \ref{sec:mms}).
\begin{open}
Does there always exist an MMS-fair allocation of a multicake with $m\geq 2$ islands among $n\geq 2$ agents with $k\geq 2$ intervals per agent?
\end{open}

\fi

\section{Envy-free Allocation}
\label{sec:envyfree}
The goal in this section is to find a $k$-feasible allocation which not only guarantees to each agent a value above some threshold, but is also \emph{envy-free}: each agent should value his/her piece at least as much as the piece of every other agent.
Formally, we want that for every two agents $i,j\in[n]$: $V_i(Z_i)\geq V_i(Z_j)$.

\renewcommand{\vbestk}[2]{U^k_{#1}(#2)} 
Define the \emph{utility} that an agent $i\in[n]$ gains from a piece $Z\subseteq C$ as
\begin{align*}
\vbestk{i}{Z} := V_i(\bestk{i}(Z)).
\end{align*}
Observe that
any general allocation can be converted to a $k$-feasible allocation by letting each agent $i$ pick the $k$ most valuable connected pieces contained in $Z_i$; 
in this case, the agent enjoys a utility of $\vbestk{i}{Z_i}$. 
Therefore, to find a $k$-feasible envy-free allocation, it is sufficient to find a \emph{general} allocation in which for every two agents $i,j\in[n]$: $\vbestk{i}{Z_i}\geq \vbestk{i}{Z_j}$.

The proofs in this section use the following technical lemma.
\begin{lemma}
\label{lem:minvalue}
Suppose a piece $Z_i$ is a union of $m_i$ intervals, and $V_i(Z_i)\geq \max(k,m_i)$. Then $\vbestk{i}{Z_i}\geq k$. In words: agent $i$ can take from $Z_i$ at most $k$ intervals with a total value of at least $k$.
\end{lemma}
\begin{proof}
If $m_i\leq k$, then $\bestk{i}(Z_i) = Z_i$ and its value is at least $k$ by assumption.

If $m_i > k$, then the value of the entire $Z_i$ is at least $m_i$, so the average value per interval is at least $1$. By the pigeonhole principle, the $k$ most valuable intervals in $Z_i$ are worth together at least $k$.
\end{proof}

\subsection{Two Agents}
\label{sub:ef-2}
To prove Theorem \ref{thm:envyfree}(a), we present problem \ref{prob:ef2}, which is an envy-free variant of problem \ref{prob:multicake-normalized} for $n=2$ agents.
The normalization is justified by arguments similar to Lemma \ref{lem:normalization}.

\begin{problem}
\caption{
\label{prob:ef2}
Envy-Free Multicake Division for $n=2$ Agents
}
\begin{algorithmic}
\REQUIRE ~\\
\begin{itemize}
\item Positive integers $m, k$ with $m\geq nk-n+1 = 2k-1$;
\item A multicake $C$ with $m$ islands; 
\item Two value-measures $V_A,V_B$ on $C$ having $V_i(C) \geq m+n-1 = m+1$ for $i\in\{A,B\}$.
\end{itemize}

\ENSURE An allocation $\mathbf{Z} = (Z_A,Z_B)$ in which:
\begin{align*}
\vbestk{A}{Z_A}\geq k; && \vbestk{A}{Z_A}\geq \vbestk{A}{Z_B}
\\
\vbestk{B}{Z_B}\geq k; && \vbestk{B}{Z_B}\geq \vbestk{B}{Z_A}
\end{align*}
\end{algorithmic}
\end{problem}

\begin{proof}[Proof of Theorem \ref{thm:envyfree}(a)]
An algorithm for solving problem \ref{prob:ef2} is presented below.

One of the agents, say Alice, evaluates all $m$ islands, orders them in increasing order of their value, and arranges them on the interval $(0,m)$ in the following way.
The first (least-valuable) island is mapped to the leftmost sub-interval $(0,1)$. The second island is mapped to the rightmost sub-interval $(m-1,m)$. The third is mapped to $(1,2)$, the fourth is mapped to $(m-2,m-1)$, and so on. 
So, the less valuable islands are mapped to the left and right ends of $(0,m)$ alternately, while the more valuable islands are mapped to the center of $(0,m)$.

Alice cuts the interval $(0,m)$ into two pieces that are equivalent in her eyes, i.e., the best $k$ intervals in the leftmost half have the same value for her as the best $k$ intervals in the rightmost half. 
Formally, Alice picks some $x_A\in(0,m)$ such that $\vbestk{A}{0,x_A} = \vbestk{A}{x_A,m}$. There exists such $x_A$ by the intermediate value theorem, since $\vbestk{A}{0,x}$ is a continuous function of $x$ that increases from $0$ towards $\vbestk{A}{C}$, and $\vbestk{A}{x,m}$ is a continuous function of $x$ that decreases from $\vbestk{A}{C}$ towards $0$. An illustration is shown below, where $m=7$ and $k=3$:
\begin{center}
\includegraphics[width=0.7\textheight]{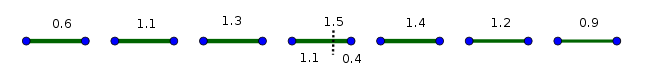}
\end{center}
The number above each island denotes its value for Alice. The dotted line denotes Alice's cut $x_A$. It divides the multicake into two parts with $\vbestk{A}{0,x_A} = \vbestk{A}{x_A,7} = 3.5$.

The other agent, say Ben, chooses the half that he prefers (the half with the larger $\vbestk{B}{\cdot}$) and takes his best $k$ intervals from it, and Alice takes her best $k$ intervals from the remaining half. This obviously yields an envy-free allocation. It remains to prove the value guarantee. 

We first prove an auxiliary claim: we prove that Alice can always make her cut in $[k-1,m-k+1]$, i.e., she never has to cut inside one of the $k-1$ leftmost islands or inside one of the $k-1$ rightmost islands. \emph{Proof of auxiliary claim.}%
\ifdefined\ACKS
\footnote{
Based on ideas in https://math.stackexchange.com/q/3045731/29780 by lulu and Gregory Nisbet.
}
\fi
~~Consider first the leftmost islands --- those mapped to $(0,k-1)$. For each such island, the next (more valuable) island in Alice's ordering is mapped to $(m-k+1,m)$. Therefore, $V_A(0,k-1)\leq V_A(m-k+1,m)$. Since both these intervals overlap less than $k$ islands, $\vbestk{A}{0,k-1}\leq \vbestk{A}{m-k+1,m}$ too. 
The assumption $m\geq 2k-1$ implies  $k-1 < k\leq  m-k+1$, so $\vbestk{A}{0,k-1}\leq \vbestk{A}{k-1,m}$. Hence Alice has a halving-point to the right of $k-1$ and can cut at $x_A\geq k-1$.
A similar consideration applies to the rightmost islands --- those mapped to $(m-k+1,m)$. For each such island, the next island in Alice's ordering is mapped to $(1,k)$. Therefore, $V_A(m-k+1,m)\leq V_A(1,k)\leq V_A(0,k)$. Since both these intervals overlap at most $k$ islands, $\vbestk{A}{m-k+1,m}\leq \vbestk{A}{0,k}$ too. 
The assumption $m\geq 2k-1$ implies that $m-k+1\geq k$, so $\vbestk{A}{m-k+1,m}\leq \vbestk{A}{0,m-k+1}$.
Hence Alice has a halving-point to the left of $m-k+1$ and can cut at $x_A\leq m-k+1$. This completes the proof of the auxiliary claim.

The auxiliary claim implies that we can \emph{require} Alice to make her cut in $[k-1,m-k+1]$, so that each half contains at least $k-1$ whole islands. 
Denote the two halves by $Z_L := (0,x_A)$ and $Z_R := (x_A,m)$ and their island counts by $m_L$ and $m_R$
with $m_L\geq k-1$ and $m_R\geq k-1$.
We now prove that, for each agent $i$, there is at least one $j\in\{L,R\}$ for which $\vbestk{i}{Z_j}\geq k$.
Since the division is envy-free, this will imply that each agent's final utility is at least $k$.

To prove the last claim we consider two cases.

\emph{Case 1}: Alice's cut is in $(k-1,m-k+1)$. Then each half contains at least $k$ sub-intervals, so $k\leq m_j$ for all $j\in\{L,R\}$.
Since Alice has made a single cut, $m_L+m_R \leq m+1$.
Recall that the valuations are normalized such that the total multicake value for each agent is at least $m+1$,
so for each agent $i$, $m_L+m_R \leq V_i(Z_L)+V_i(Z_R)$. 
Hence there exists some index $j\in\{L,R\}$ such that $m_j\leq V_i(Z_j)$.
Now $k\leq m_j$ implies $V_i(Z_j)\geq \max(k,m_j)$, and by Lemma \ref{lem:minvalue}, $\vbestk{i}{Z_j}\geq k$.

\emph{Case 2:} The edge case in which Alice's cut is at $x_A = k-1$ (the other edge case $x_A = m-k+1$ is analogous). Then $m_L = k-1$ and $m_R = m-k+1 \geq k$.
Since Alice's cut is at an integer point, $m_L + m_R = m$,
so $(m_L+1)+m_R = m+1 \leq  V_i(Z_L)+V_i(Z_R)$.
Now for each agent $i$,
either $(m_L+1)\leq V_i(Z_L)$, or $m_R\leq V_i(Z_R)$.
In the former case, 
since $Z_L$ contains less than $k$ islands,
$\vbestk{i}{Z_L} = V_i(Z_L) \geq m_L+1 = k$.
In the latter case,
since $m_R\geq k$, 
$V_i(Z_R)\geq \max(m_R,k)$, so 
by Lemma \ref{lem:minvalue}, $\vbestk{i}{Z_R}\geq k$.
\end{proof}

\subsection{Many agents}
\label{sub:ef-n}
To prove Theorem \ref{thm:envyfree}(b),
we present problem \ref{prob:efn}, which is an envy-free variant of problem \ref{prob:multicake-normalized} differing only in the normalization.
Note that the lower bound of the theorem, namely $\min\left({1\over n},{k\over nk+m-k}\right)$, equals 
${k\over nk+\max(k,m)-k}$ by elementary arithmetics.

\begin{problem}
\caption{
\label{prob:efn}
Envy-Free Multicake Division for Many Agents.
}
\begin{algorithmic}
\REQUIRE ~\\
\begin{itemize}
\item Positive integers $m, n, k$;
\item A multicake $C$ with $m$ islands; 
\item $n$ value-measures $V_1,\ldots,V_n$ on $C$ having $V_i(C) \geq nk + \max(k,m) - k$ for all $i\in[n]$.
\end{itemize}

\ENSURE An allocation $\mathbf{Z} = Z_1,\ldots,Z_n$ in which, for each agent $i\in[n]$:
\begin{align*}
\vbestk{i}{Z_i}\geq k;
&&
\vbestk{i}{Z_i}\geq \vbestk{i}{Z_j}\text{~for all $j\in[n]$.}
\end{align*}
\end{algorithmic}
\end{problem}

\begin{proof}[Proof of Theorem \ref{thm:envyfree}(b)]
An algorithm for solving problem \ref{prob:efn} is presented below.

Arrange the islands in an arbitrary order into the $m$ intervals $(0,1),\ldots,(m-1,m)$.

Consider any partition of the interval $(0,m)$ into $n$ intervals $Z_1,\ldots, Z_n$. In each such partition, if an agent is asked to choose a preferred piece, he selects the piece $Z_j$ with the highest utility, i.e., $\arg\max_{j\in[n]} \vbestk{i}{Z_j}$.
This selection satisfies two properties: (1) An agent always weakly prefers a nonempty piece over an empty piece; (2) The set of partitions in which an agent prefers the piece with index $j$ is a closed subset of the space of all partitions. This is because $\vbestk{i}{[x,y]}$ is a continuous function of $x$ and $y$.
In their classic papers,
\citet{Stromquist1980How} and \citet{Su1999Rental} proved that, whenever the preferences of $n$ agents satisfy these two properties, there exists a connected partition of $(0,m)$, in which each agent prefers a different piece, i.e., the partition is envy-free.
Let $(Z_1,\ldots,Z_n)$ be such a partition.

It remains to prove that, in this envy-free partition,
$\vbestk{i}{Z_i}\geq k$ for all $i$.
In words: each agent $i$ can find in his/her share $Z_i$ at most $k$ connected pieces with a total value of at least $k$.
Since the allocation is envy-free, it is sufficient to prove that, for all $i\in[n]$, there exists {some} $j\in[n]$ with $\vbestk{i}{Z_j}\geq k$.

The proof focuses on a specific agent $i$.  We first discard all the parts $Z_j$ for which $V_i(Z_j)<k$. 
Suppose $n-d$ such parts are discarded (for some $d\range{0}{n}$). 
Recall that the total multicake value is normalized to $nk + \max(k,m) - k$; therefore the total value of the remaining $d$ parts is strictly larger than $(d-1)k + \max(k,m)$. This expression is at least 0, so the remaining value is strictly larger than 0, so at least one part remains, so in fact $d\geq 1$.

Suppose w.l.o.g. that the remaining parts are $Z_1,\ldots,Z_d$, and each part $Z_j$ contains $m_j$ sub-intervals. 
The total number of sub-intervals in $Z_1,\ldots,Z_n$ is at most $m+n-1$ (since there were $m$ islands in the input instance, and the algorithm made $n-1$ cuts). After removing $n-d$ parts, each of which contains at least one interval, the total number of sub-intervals is at most $d+m-1$:
\begin{align}
\label{eq:dm1a}
\sum_{j=1}^{d} m_j ~\leq~ d+m-1.
\end{align}
On the other hand, the total value of the remaining parts is:
\begin{align}
\label{eq:dm1b}
\sum_{j=1}^d V_i(Z_j) ~>~ (d-1)\cdot k + \max(k,m) \geq (d-1)\cdot 1 + m = d+m-1
\end{align}
Comparing \eqref{eq:dm1a} and \eqref{eq:dm1b} implies that $\sum_{j=1}^d V_i(Z_j) \geq \sum_{j=1}^d m_j$, so
for at least one $j\range{1}{d}$:
\begin{align*}
V_i(Z_j)\geq m_j.
\end{align*}
Moreover, since we started by discarding the parts whose value is less than $k$, in fact: $V_i(Z_j)\geq \max(k,m_j)$. Hence, by Lemma \ref{lem:minvalue}, $\vbestk{i}{Z_j}\geq k$. Agent $i$ does not envy agent $j$, so
$\vbestk{i}{Z_i}\geq k$ too.
\end{proof}

\subsection{Open questions}
When either $n=2$ or $k=1$ (or both), the results presented in this section are tight, and they imply that $r_E(m,n,k)=r(m,n,k)$.
However, for $n\geq 3$ and $k\geq 2$ the results are not tight.
Intuitively, the algorithm of \S\ref{sub:ef-n} does not use all the freedom allowed by the problem: it begins by ordering the islands arbitrarily on a line, and then gives each agent a contiguous subset of that line.
The following example shows that this approach cannot guarantee a larger value.

Suppose all agents have the same value-measure: they value the $m-1$ leftmost islands $(0,1),\ldots, (m-2,m-1)$ at $1$, and the rightmost island $(m-1,m)$ at $nk-k+1$, so the total multicake value is $nk+m-k$. Now there are three cases:
\begin{enumerate}
\item At least one agent gets all his/her $k$ pieces in $k$ of the islands $1,\ldots,m-1$. This agent receives a value of at most $k$.
\item All $n$ agents get a piece of island $m$, and the leftmost agent (who gets the leftmost interval of $(0,m)$) gets a value of at most $1$ from island $m$. This leftmost agent can get at most $k-1$ other islands and remains with a total value of at most $k$.
\item All $n$ agents get a piece of island $m$, and the leftmost agent gets a value of more than $1$ from island $m$. 
Then, the remaining value in island $m$ is less than $nk-k = (n-1)k$. Therefore, at least one of the $n-1$ remaining agents gets a value of less than $k$.
\end{enumerate}
To overcome this impossibility, it may be required to re-order the islands based on their value as in \S\ref{sub:ef-2}. It is not clear how to do this while keeping the envy-freeness guarantee.
In theory, there could be a separation between the attainable value with and without envy.

\begin{open}
With $n\geq 3$ agents and $k\geq 2$ pieces per agent, what value can be guaranteed to all agents in an envy-free allocation?
Is it smaller than the value-guarantee in an allocation that may have envy?
\end{open}

The proof of Theorem \ref{thm:envyfree}(b) is existential --- it uses results on the existence of a connected envy-free allocation of an interval. 
It is known that, even in the classic setting where $m=1$ and $k=1$, such allocations cannot be found by a finite protocol when $n\geq 3$. \citep{Stromquist2008Envyfree}. However, approximate envy-free allocations can be computed efficiently. Recently, \citet{arunachaleswaran2019fair} presented a multiplicative approximation algorithm and \citet{goldberg2020contiguous} presented some additive approximation algorithms. However, it is not immediately clear that these algorithms can be applied in our setting, since they assume that the agents' valuations are additive, while the functions $\vbestk{i}{\cdot}$ are not additive.%
\footnote{
The $\vbestk{i}{\cdot}$ are not even subadditive. For example, suppose $k=1$ and $Z_1 = [0,0.1]\cup[0.2,0.3]$ and $Z_2 = [0.1,0.2]$. Then $\vbestk{i}{Z_1\cup Z_2} = V_i(Z_1\cup Z_2)=V_i(Z_1)+V_i(Z_2)$, which may be larger than $\vbestk{i}{Z_1}+\vbestk{i}{Z_2}$.
} 
\begin{open}
Can approximate envy-free multicake allocations can be found by an efficient algorithm?
\end{open}

Finally, let us mention an open problem that is interesting even in the classic setting of a single interval ($m=1$).
The impossibility result of \citet{Stromquist2008Envyfree} is proved only for $k=1$.
Allowing $k = 2$ intervals per agent, the famous Selfridge-Conway protocol finds an envy-free cake-allocation in bounded time, but only for $n = 3$ agents. In contrast, the recent protocols by \citet{Aziz2015Discrete} and \citet{amanatidis2018improved} find 
an envy-free cake-allocation for $n = 4$ agents, but may allocate each agent much more than two pieces.
\begin{open}
Does there exist a finite protocol that finds an envy-free allocation of a single interval among $n\geq 4$ agents giving at most $k=2$ intervals to each agent?
\end{open}

\section{Dividing a Rectilinear Polygon}
\label{sec:rectilinear}

\begin{figure}
\begin{center}
\includegraphics[width=.25\textwidth]{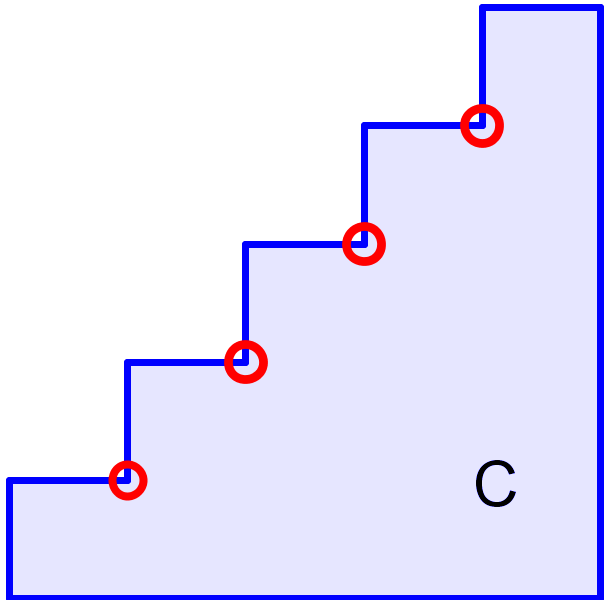}
\hskip 5cm	
\includegraphics[width=.25\textwidth]{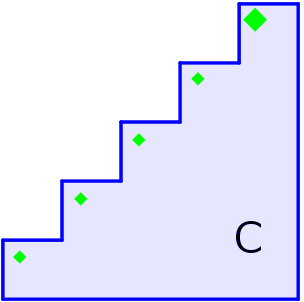}
\end{center}
\caption{
\label{fig:rectilinear}
\textbf{Left}: a rectilinear polygon. Circles denote reflex vertices.
\textbf{Right}: the value function of Proposition \ref{prop:rectilinear-negative}.
}
\end{figure}
This section shows an application of the multicake division algorithm for fairly dividing a two-dimensional land-estate.
A natural requirement in land division settings is to give each agent a rectangular piece \citep{SegalHalevi2020EnvyFree}. 
We assume that the cake $C$ is a \emph{rectilinear polygon} --- a polygon whose angles are $90^\circ$ or $270^\circ$. 
Rectilinearity is a common assumption in practical partition problems \citep{Keil2000Polygon}.  The complexity of a rectilinear polygon is characterized by the number of its \emph{reflex vertices} --- vertices with a $270^\circ$ angle. We denote this number by $T$. A rectangle --- the simplest rectilinear polygon --- has $T=0$. The polygon in Figure \ref{fig:rectilinear} has $T=4$ reflex vertices.

\begin{theorem}
	\label{thm:rectilinear}
	It is possible to divide a rectilinear polygon with $T$ reflex vertices among $n$ agents giving each agent a rectangle with value at least 
$
	1 / (n+T)
$
	of the total polygon value.
\end{theorem}
\begin{proof}
	\citet{Keil2000Polygon,Eppstein2010GraphTheoretic}
present efficient algorithms for partitioning a rectilinear polygon into a minimal number of rectangles. A rectilinear polygon with $T$ reflex vertices can be partitioned in time $O(\text{poly}(T))$ into at most $T+1$ rectangles (this number is tight when the polygon vertices are in general position). 
Using such algorithms, partition the polygon $C$ into at most $T+1$ rectangular ``islands'' and then apply the theorems for multicake division with $m=T+1$. By setting $k=1$ in Theorem \ref{thm:prop}, each agent $i$ gets a single rectangle worth at least $V_i(C)/(n+T)$.
\end{proof}

If each agent is willing to receive up to $k$ rectangles, then the value-guarantee (as a fraction of the total multicake value) increases accordingly to $\min~(1/n,~k/(n+T))$.

Theorem \ref{thm:rectilinear} is tight in the following sense: 
\begin{proposition}
\label{prop:rectilinear-negative}
for every integer $T\geq 0$, there exists a rectilinear polygon with $T$ reflex vertices, in which it is impossible to guarantee every agent a rectangle with a value of more than $1/(n+T)$ the total polygon value. 
\end{proposition}
\begin{proof}
Consider a staircase-shaped polygon with $T+1$ stairs as in Figure \ref{fig:rectilinear}/Right. 
All agents have the same value-measure, which is concentrated in the diamond-shapes: the top diamond is worth $n$ and each of the other diamonds is worth 1. So for all agents, the total polygon value is $n+T$.
Any rectangle in $C$ can touch at most a single diamond. There are two cases: 
\begin{itemize}
\item At least one agent touches one of the $T$ bottom diamonds. Then, the value of that agent is at most 1. 
\item All $n$ agents touch the top diamond. Then, their total value is $n$ and at least one of them must receive a value of at most 1.
\end{itemize}
So, either some or all of the agents receive a value of at most $V_i(C)/(n+T)$.
\end{proof}

\subsection{Open question}
While Theorem \ref{thm:rectilinear} is tight for a worst-case polygon, it is not necessarily tight for specific polygons. For example, consider the following polygon:
\begin{center}
	\includegraphics[width=.3\textwidth]{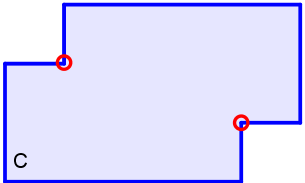}
\end{center}
It has $T=2$ reflex vertices, can be partitioned into $T+1=3$ disjoint rectangles, and cannot be partitioned into less disjoint rectangles. Hence, for a single agent, Theorem \ref{thm:rectilinear} guarantees $1/(1+2) = 1/3$ of the total value. However, since the polygon can be covered by $2$ overlapping rectangles, at least one of them is worth at least $1/2$ of the total value, so the actual value-guarantee for a single agent is $1/2$. 
In general, for $n=1$ agent, Theorem \ref{thm:rectilinear} guarantees $V_i(C)$ divided by the \emph{partition-number} of $C$ (the number of pairwise-disjoint rectangles whose union equals $C$), while it is possible to guarantee 
$V_i(C)$ divided by the \emph{cover-number} of $C$ (the number of possibly-overlapping rectangles whose union equals $C$).

\begin{open}
Is it possible to divide a rectilinear polygon among $n\geq 2$ agents,
with a value-guarantee that depends on its rectangle cover number rather than on $T$?
\end{open}

\section*{Acknowledgments}
I am grateful to Elad Aigner-Horev, Chris Culter\footnote{in \url{http://math.stackexchange.com/q/461675}}, Zur Luria, Yuval Filmus, Gregory Nisbet and lulu for their helpful ideas,
and three anonymous reviewers of Discrete Applied Mathematics for their many constructive comments.
The research was partly funded by the Doctoral Fellowships of Excellence Program at Bar-Ilan University, the Mordecai and Monique Katz Graduate Fellowship Program, and the Israel Science Foundation grant 1083/13.

\appendix

\section*{APPENDIX}
\section{Finding an envy-free matching in a bipartite graph}
\label{sec:efm}
This appendix presents a simplified proof of a corollary that is proved by \citet{aigner2019envy}.

For the algorithm below, recall that, given a matching $M$ in a graph $G$, an \emph{$M$-alternating path} is a path in $G$ whose first edge is not in $M$, second edge is in $M$, third edge is not in $M$, and so on.

\begin{algorithm}
\caption{
\label{alg:efm}
Finding a nonempty envy-free matching (cf. Definition \ref{def:efm}).
}
\begin{algorithmic}[1]
\REQUIRE A bipartite graph $G = (X\cup Y,E)$ satisfying $|N_G(X)| \geq |X| \geq 1$.

\ENSURE An envy-free matching between nonempty subsets $X_L\subseteq X$ and $Y_L\subseteq Y$.

\dottedline{}

\STATE 
\label{step:maxmatch}
Let $M$ be a maximum-cardinality matching in $G$.

\STATE 
\label{step:unmatched}
Let $X_0$ be the set of vertices of $X$ unmatched by $M$.

\STATE 
\label{step:reachable}
Let $X_S, Y_S$ be the subsets of $X, Y$ respectively reachable from $X_0$ by $M$-alternating paths.

\RETURN the subset of $M$ adjacent to $X_L := X\setminus X_S$.
\end{algorithmic}
\end{algorithm}

\begin{theorem}
\label{thm:efm}
Algorithm \ref{alg:efm}:

(a) Runs in polynomial time;

(b) Returns an envy-free matching;

(c) Returns a nonempty matching whenever its input conditions are satisfied.
\end{theorem}
\begin{proof}
(a) In step \ref{step:maxmatch}, a maximum-cardinality matching can be find in polynomial time, e.g. by the Hopcroft-Karp algorithm.
Step \ref{step:unmatched} is obviously linear.
In step \ref{step:reachable}, the sets of reachable vertices can be found in polynomial time using breadth-first search starting at $X_0$.

(b) All vertices of $X_L$ are matched by $M$, since by definition the unmatched vertices are in $X_0$ which is contained in $X_S$. Hence the vertices of $X_L$ are not envious.
Suppose by contradiction that a vertex $x_0\in X_S$ is envious. By definition of envy-free matching, this $x_0$ is adjacent to some vertex $y\in Y$,
which is matched by $M$ to some other vertex $x_1\in X_L$.
But this means that $x_l$ is reachable from $X_0$ by an $M$-alternating path (through $x_0$ and $y$), so $x_1\in X_S$. This contradicts the definition of $X_L$.

(c) If $X_0$ is empty, then $X_S$ is empty so $X_L = X$. By the input assumption, $|X|\geq 1$.

Suppose $X_0$ is not empty. By construction of $X_S$ and $Y_S$, every vertex in $Y_S$ is matched by $M$ to some vertex in $X_S$
(if a vertex $y\in Y_S$ were not matched by $M$, then we could trace an $M$-augmenting path from $y$ to $x_0$, contradicting the maximality of $M$). However, the vertices of $X_0$ are by definition not matched by $M$. Therefore, $|X_S| > |Y_S|$. Now, by construction of $X_S$ and $Y_S$, all neighbors of $X_S$ are in $Y_S$: $N_G(X_S)\subseteq Y_S$. Therefore, 
$|N_G(X_S)| < |X_S|$. By the input conditions, $|N_G(X)| \geq |X|$. Therefore $X_S\neq X$, so $X_L$ is not empty.
\end{proof}

\section{Maximin-share fairness}
\label{sec:mms}
Given a multicake $C$ and integers $k\geq 1$ and $n\geq 1$, the \emph{maximin share (MMS)} of agent $i$, denoted $\mmsk{i}(C)$, is the largest value that $i$ can get by
constructing a $k$-feasible $n$-partition of (a subset of) $C$, and getting the least valuable part in that partition.
An allocation is \emph{MMS-fair} if the value allocated to each agent $i$ is at least $\mmsk{i}(C)$.

As in the previous algorithms, we present a normalized problem.
\begin{problem}
\caption{
\label{prob:mms-normalized}
MMS-fair Multicake Division Problem
}
\begin{algorithmic}
\REQUIRE ~\\
\begin{itemize}
\item Positive integers $m, n, k$; 
\item A multicake $C$ with $m$ islands; 
\item $n$ value-measures $V_1,\ldots,V_n$ having  $\mmsk{i}(C) \geq k$ for all $i\in[n]$.
\end{itemize}

\ENSURE A $k$-feasible allocation $\mathbf{Z} = Z_1,\ldots,Z_n$ in which for all $i\in[n]$:
$
V_i(Z_i)\geq k.
$
\end{algorithmic}
\end{problem}

The normalization implies that, for each agent $i\in[n]$, there exist $n$ pieces $Y_{i,1},\ldots,Y_{i,n}$, each of which contains at most $k$ intervals, such that, for all $j\in[n]$, $V_i(Y_{i,j})\geq k$. 


\begin{maintheorem}
\label{thm:mms}
When $k=1$, for any $n\geq 1$ and $m\geq 1$, Algorithm \ref{alg:main} solves Problem \ref{prob:mms-normalized}.
\end{maintheorem}
\begin{proof}
By induction on $n$. If $n=1$ then the algorithm terminates at step 1. 
The MMS allocation of the single remaining agent is made of a single contiguous interval with a value of at least $1$, so the most valuable island is worth for this agent at least $1$.

Suppose now that $n>1$ and that the claim holds whenever there are less than $n$ agents.
At step 2 the algorithm finds a partial allocation. When $k=1$, such partial allocation can always be found by a mark auction (Algorithm \ref{alg:mark}),
with $A = \emptyset$ as a barren subset of $k-1$ islands,
and $B = \bestk{1}(C) = $ the island most valuable to agent \#1.
The normalization guarantees that island $B$ is worth at least $1$ for agent \#1.
The mark auction does not depend in any way on the normalization, so it still (by Lemma \ref{lem:mark}) returns a partial allocation.
In this allocation, a single agent, say agent $n$, receives a single interval $Z_n$, which is at the leftmost end of the island $B$. This $Z_n$ is worth exactly $1$ for agent $n$, so the output condition is satisfied for $n$.

Let $C' := C\setminus Z_n$ be the remaining multicake. 
We have to show that, for the remaining $n-1$ agents, the input requirements of Problem \ref{prob:mms-normalized} are still satisfied, i.e., for every agent $i \in [n-1]$, $\mmskm{i}(C') \geq 1$.

By the input conditions, before the piece $Z_n$ was allocated, there existed $n$ intervals $Y_{i,1},\ldots,Y_{i,n}$ in $C$, such that, for all $j\in[n]$, $V_i(Y_{i,j})\geq 1$.
Now, $Z_n$ is the leftmost part of the island $B$, 
and $V_i(Z_n)\leq 1$, 
so $Z_n$ overlaps at most one of these $n$ intervals, say $Y_{i,n}$.
This means that, for all $j\in[n-1]$, the intervals $Y_{i,n-1}$ do not overlap $Z_n$, and thus they are all subsets of $C'$.
Therefore, the $n-1$ intervals $Y_{i,1},\ldots,Y_{i,n-1}$ show that $\mmskm{i}(C') \geq 1$.
By the induction assumption, the recursive call solves problem \ref{prob:mms-normalized} for the remaining agents.
\end{proof}

Unfortunately, 
when $k=2$, Algorithm \ref{alg:main} may fail to find an MMS-fair allocation. 
As an example, consider a multicake with $m=4$ islands, denoted $C_1,C_2,C_3,C_4$, and $n=2$ agents who value the islands at $1.5, 1.5, 0.5, 0.5$ respectively. 
The MMS of both agents is $2$, by the partition $(C_1\cup C_3, C_2\cup C_4)$.
All islands are barren, so Algorithm \ref{alg:main} starts by running a mark auction, with e.g. $A = C_1$ and $B = C_2$. Both agents cut $C_2$ at the same point; one of them receives $C_1$ and the leftmost part of $C_2$, which is worth $0.5$. 
The remaining cake now contains 
the remainder of $C_2$, which is worth $1$,
and two more islands with values $0.5, 0.5$. 
Now the remaining agent cannot get two intervals with a value of $2$.

In this particular case, an MMS allocation obviously exists, but in general,
when 
 $n\geq 2$ and $k\geq 2$,  the existence of an MMS-fair multicake allocation remains open.

 \bibliographystyle{ACM-Reference-Format}
 \bibliography{../erelsegal-halevi}

\end{document}